\documentclass[10pt]{article}

\usepackage{amsmath}
\usepackage{amsfonts} 
\usepackage{amssymb,amsthm}
\usepackage{color}
\definecolor{red}{rgb}{1,0,0}
\definecolor{gray}{gray}{0.5}
\definecolor{darkblue}{rgb}{0,0,0}
\usepackage{hyperref}

\usepackage{graphicx}

\usepackage[english]{babel}   

\usepackage{pb-diagram}  
\usepackage[cmtip,arrow]{xy}
\usepackage{pb-xy} 
\usepackage{enumerate}   
\usepackage{bbold} 
\usepackage[utf8]{inputenc} 
\usepackage{stmaryrd}

\newtheorem{de}{Definition}[section]
\newtheorem{theo}[de]{Theorem}
\newtheorem{prop}[de]{Proposition}

\newtheorem{lem}[de]{Lemma}

\newcounter{qcounter}

\newenvironment{noindlist}
 {\begin{list}{\textbf{(\arabic{section}.\arabic{qcounter})} }{\usecounter{qcounter} \leftmargin=0em 
\labelwidth=0mm \labelsep=0cm  }}
 {\end{list}}

\newcommand{\bb}[1][R]{\mathbb{#1}}
\newcommand{\sbb}[1][1|1]{\mathbb{R}^{#1}}

\newcommand{\bbc}[1][C]{\mathbb{#1}}
\newcommand{\sbbc}[1][1|1]{\mathbb{C}^{#1}}

\newcommand{\sbbpq}{\mathbb{R}^{p|q}}
\newcommand{\sbbmn}{\mathbb{R}^{m|n}}
\newcommand{\ca}[1][M]{{\mathcal{#1}}} 
\newcommand{\can}[1][N]{{\mathcal{#1}}}

\newcommand{\om}[1][M]{{\ca[O]_{\ca[#1]}}}

\newcommand{\tm}[1][M]{\ca[T]_{\ca[#1]}}

\newcommand{\pa}[1][q_i]{\partial_{#1}}
\newcommand{\sig}[1][|q_i|]{(-1)^{#1}}

\newcommand{\cat}[1][M]{\ca[#1]=(#1,\om[#1])}

\newcommand{\dis}{\displaystyle}
  \newcommand{\sk}[1][0|k]{\mathbb{R}^{#1}}

 \newcommand{\cg}[1][G]{{\mathcal{#1}}} 

\begin{document}

\title{Integration of vector fields on smooth\\
and holomorphic supermanifolds}
\author{
{\bf St\'ephane Garnier}\\
     {\small  Fakult\"at f\"ur Mathematik}\\
             {\small Ruhr-Universit\"at Bochum}\\
                 {\small D-44780 Bochum, Germany}\\
                     {\small \tt{Stephane.Garnier@rub.de}\vspace*{3mm}}\\
 {\small and\vspace*{3mm}}\\
{ \bf Tilmann Wurzbacher\footnote{On leave of absence; present address: Fakult\"at f\"ur Mathematik, Ruhr-Universit\"at
Bochum, D-44780 Bochum, Germany; e-mail: \tt{Tilmann.Wurzbacher@rub.de}}}\\  
   {\small LMAM-UMR 7122}\\
 {\small Universit\'e de Lorraine et }\\
 {\small C.N.R.S.}\\
 {\small F-57045 Metz, France}\\
{\small  \tt{wurzbacher@math.univ-metz.fr}}
 }

\date{\today}

\maketitle

\centerline{ {\bf Keywords: } Supermanifolds; Vector fields; Flows; Group actions}

\vskip0.3cm

\centerline{  {\bf MSC2010:} Primary 58A50, 37C10; Secondary 57S20, 32C11
}

\begin{abstract}
We give a new and self-contained proof of the existence and unicity of the flow for an arbitrary (not necessarily 
homogeneous) smooth vector field on a real supermanifold, and extend these results to the case of 
holomorphic vector fields on complex
supermanifolds. Furthermore we discuss local actions associated to super vector fields, and give several examples
and applications, as, e.g., the construction of an exponential morphism for an arbitrary finite-dimensional Lie
supergroup.
\end{abstract}

\section{Introduction}
The natural problem of integrating vector fields to obtain appropriate ``flow maps'' on supermanifolds is
considered in many articles and monographs (compare, e.g., \cite{Shander}, \cite{Bruzzo},
\cite{Tuynman}, \cite{Rogers} and \cite{dumi}) but a
``general answer'' was to our knowledge only given in the work of J. Monterde and co-workers (see
\cite{Monterde} and
\cite{Mo-SV2}). Let us consider a supermanifold $\cat$ together with a vector
field $X$ in $\tm(M)$, and an initial condition $\phi$ in $\text{Mor}(\ca[S],\ca)$,
where $\ca[S]=(S,\om[S])$ is an arbitrary supermanifold and $\text{Mor}(\ca[S],\ca)$ denotes the set of 
morphisms from $\ca[S]$ to $\ca$. The case of
classical, ungraded, manifolds leads one to consider the following question: does there
exist a ``flow map'' $F$ defined on an open sub supermanifold
$\ca[V]\subset\sbb\times\ca[S]$ and having values in $\ca$ and
an appropriate derivation on $\sbb$,
$D=\pa[t]+\pa[\tau]+\tau(a\pa[t]+b\pa[\tau])$, where $\pa[t]=\frac{\partial}{\partial t}$,
$\pa[\tau]=\frac{\partial}{\partial \tau}$ and $a$, $b$ 
are real numbers, such that the following equations are fulfilled
\begin{equation}
 \begin{array}{lcl}
D\circ F^* & = & F^*\circ X\:\\
F\circ\text{inj}^{\ca[V]}_{\{0\}\times\ca[S]} & = & \phi\:.\label{eqflowintro}
\end{array}
\end{equation}
Of course, $\ca[V]$ should be a ``flow domain'', i.e. an open sub supermanifold
of $\sbb\times\ca[S]$ such that $\{0\}\times S$ is contained in the body $V$ of $\ca[V]$
and for $x$ in $S$, the set $I_x\subset\bb$ defined by
$I_x\times\{x\}=(\bb\times\{x\})\cap V$ is an open interval. Furthermore $\text{inj}^{\ca[V]}_{\{0\}\times\ca[S]}$
denotes the natural injection morphism of the closed sub supermanifold $\{0\}\times\ca[S]$ of $\ca[V]$ into $\ca[V]$.
Of course, we could concentrate on the case $\ca[S]=\ca$ and
$\phi=\text{id}_{\ca}$, but it will be useful for our later arguments to state
all results in this (formally) more general setting.\\

Though for homogeneous vector fields ($X=X_0$ or $X=X_1$) system
($\ref{eqflowintro}$) does always have a solution, in the general case ($X=X_0+X_1$ with
$X_0\neq 0$ and $X_1\neq 0$) the system is overdetermined. A simple example of an inhomogeneous vector
field such that (\ref{eqflowintro}) is not solvable is given by $X=X_0+X_1=\left( \frac{\partial}{\partial
x}+\xi\frac{\partial}{\partial \xi} \right)+\left( \frac{\partial}{\partial \xi}+\xi\frac{\partial}{\partial x} \right)$
on $\ca=\sbb$. The crucial novelty of \cite{Mo-SV2} is to
consider instead of (\ref{eqflowintro}) the following modified, weakened, problem 
\begin{equation}
 \begin{array}{lcl}
(\text{inj}_{\bb}^{\sbb})^*\circ D\circ F^* & = &
(\text{inj}_{\bb}^{\sbb})^*\circ F^*\circ X\: \\
F\circ\text{inj}^{\ca[V]}_{\{0\}\times\ca[S]} & = & \phi\:,\label{eqflowintro2}
\end{array} 
\end{equation}
where $\text{inj}_{\bb}^{\sbb}=\text{inj}_{\bb\times \ca[S]}^{\sbb\times\ca[S]}$ is again the natural injection (and
where the above more general derivation $D$
could be replaced by $\pa[t]+\pa[\tau]$ 
since $(\text{inj}_{\bb}^{\sbb})^*$ annihilates germs of superfunctions of the type
$\tau\cdot f, f\in\ca[O]_{\bb\times\ca[S]})$.\\

In \cite{Mo-SV2} (making indispensable use of \cite{Monterde}) it is shown that
in the smooth case (\ref{eqflowintro2}) has a unique maximal solution $F$, defined on the
flow domain $\ca[V]=(V,\ca[O]_{\sbb\times\ca[S]}|_{V})$, where $V\subset \bb\times S$
is the maximal flow domain for the flow of the reduced vector field
$\widetilde{X}=\widetilde{X_0}$ on $M$ with initial condition
$\widetilde{\phi}$. Since the results of \cite{Monterde} are obtained by the use
of a Batchelor model for $\ca$, i.e. a real vector bundle $E\to M$ such that
$\ca\cong(M,\Gamma_{\Lambda E^*}^\infty)$, and a connection on $E$, we follow
here another road, closer to the classical, ungraded, case and also
applicable in the case of complex supermanifolds and holomorphic vector fields.
\\

Our new method of integrating smooth vector fields on a supermanifold in Section 2 
consists in first locally solving a finite hierarchy of ordinary differential equations,
and is here partly  inspired by the approach of
\cite{dumi}, where the case of homogeneous super vector fields on compact
supermanifolds is treated. We then show existence and unicity of solutions of
(\ref{eqflowintro2}) on smooth supermanifolds and easily deduce the results of
\cite{Mo-SV2} from our Lemmata \ref{theolocal} and \ref{theoglob}.\\

 A second beautiful result of \cite{Mo-SV2} (more precisely, Theorem 3.6 of that
reference) concerns the question if the flow $F$ solving (\ref{eqflowintro2})
fulfills ``flow equations'', as in the ungraded case. Hereby, we mean the
existence of a Lie supergroup structure on $\sbb$ such that $F$ is a local
action of $\sbb$ on $\ca$ (in case $\ca[S]=\ca,\phi=\text{id}_{\ca}$). Again,
the answer is a little bit unexpected: in general, given $X$ and its flow
$F:\sbb\times\ca\supset\ca[V]\to\ca$, there is no Lie supergroup structure on
$\sbb$ such that $F$ is a local $\sbb$-action (with regard to this structure).
The condition for the existence of such a structure on $\sbb$ is equivalent to
the condition that (\ref{eqflowintro2}) holds without the post-composition with
$(\text{inj}_{\bb}^{\sbb})^*$, i.e. the overdetermined system
(\ref{eqflowintro}) is solvable. Furthermore,
both conditions cited are equivalent to the condition that $\bb X_0\oplus\bb
X_1$ is a sub Lie superalgebra of $\tm(M)$, the Lie superalgebra of all vector
fields on $\ca$.\\

After discussing Lie supergroup structures and right invariant vector fields on  $\sbb$, as well
as local Lie group actions in the category of supermanifolds in general, we show in Section 3 the 
equivalence of the above three conditions, already given in \cite{Mo-SV2}. We include our proof
here notably  in order to be able to apply it in the holomorphic case in Section 5 (see below) 
by simply indicating how to adapt it to this context. Let us nevertheless observe that our result is slightly 
more general since we do not need to ask for any  normalization  of the supercommutators between $X_1$ and 
$X_0$ resp. $X_1$, thus giving the criterion some extra flexibility  in applications. \\

In Section 4, we give several examples of vector fields on supermanifolds, homogeneous and inhomogeneous, and explain
their integration to flows. Notably, we construct an exponential morphism for an arbitrary finite-dimensional Lie
supergroup, via a canonically defined vector field and its flow. We comment here also on the integration of what 
are usually called ``(infinitesimal) supersymmetries" in physics, i.e., purely odd vector fields having non-vanishing 
self-commutators. \\

Finally, in Section 5 we adapt our method to obtain flows of vector fields 
(compare Section 2 and notably Lemma \ref{theolocal}) to the case of holomorphic vector
fields on holomorphic supermanifolds. To avoid monodromy problems one has, of course, to take care of the 
topology of the flow domains, and maximal flow domains are -as already in the ungraded holomorphic case- 
no more unique. Otherwise the analogues of all results in Section 2 and 3 continue to hold in the 
holomorphic setting.\\

Throughout the whole article we will work in the ringed space-approach to supermanifolds (see, e.g.,
\cite{Kostant}, \cite{Leites}, \cite{Manin} and \cite{Schmitt} for detailed accounts of this approach). Given two
supermanifolds $\cat$ and $\cat[N]$, a ``morphism'' $\phi=(\widetilde{\phi},\phi^*):\ca\to\can$ is thus given by a
continuous map $\widetilde{\phi}:M\to N$ between the ``bodies'' of the two supermanifolds and a sheaf homomophism
$\phi^*:\om[N]\to\widetilde{\phi}_{*}\om$. The topological space $M$ comes canonically with a sheaf
$\ca[C]^\infty_M=\om/\ca[J]$, where $\ca[J]$ is the ideal sheaf generated by the germs of odd superfunctions, such that
$(M,\ca[C]^\infty_M)$ is a smooth real manifold. Then $\widetilde{\phi}$ is a smooth map from $(M,\ca[C]^\infty_M)$ to 
$(N,\ca[C]^\infty_N)$. Let us recall that a (super) vector field on $\cat$ is, by definition, an element of the Lie
superalgebra $\tm(M)=(\text{Der}_{\bb}(\om))(M)$ and that $X$ always induces a smooth vector field $\widetilde{X}$ on
$(M,\ca[C]^\infty_M)$. For $p$ in $M$ and $f+\ca[J]_{p}\in(\ca[C]^\infty_M)_p=(\om/\ca[J])_p$ one defines
$\widetilde{X}_p(f+\ca[J]_{p})=X_0(f)(p)$, where $X_0$ is the even part of $X$ and for $g\in (\om)_p$, $g(p)\in\bb$ is
the value of $g$ in the point $p$ of $M$.

\section{Flow of a vector field on a real supermanifold}

In this section we give our main result on the integration of general (i.e. not necessary homogeneous) vector fields by
a new method, avoiding auxiliary
choices of Batchelor models and connections, as in \cite{Monterde}. Our more direct approach is inspired, e.g., by
\cite{dumi}, where the
case of homogeneous
vector fields on compact manifolds is treated, and it can be adapted to the holomorphic case (see Section 4).\\ 

 For the sake of readability we will often use the following shorthand:
if $\ca[P]$ is a supermanifold, we write $\text{inj }_{\bb}^{\sbb}$ for 
$\text{inj }_{\bb\times\ca[P]}^{\sbb\times\ca[P]}$. Furthermore, the canonical coordinates of $\sbb$ will be denoted by
$t$ and $\tau$, with ensueing vector fields $\pa[t]=\frac{\partial }{\partial t}$ and $\pa[\tau]=\frac{\partial
}{\partial \tau}$.
\begin{lem}\label{theolocal}
Let $\ca[U]\subset\sbb[m|n]$ and $\ca[W]\subset\sbb[p|q]$ be superdomains,
$X\in\tm[W](W)$ be a super vector field on $\ca[W]$ (not necessarily
homogeneous) and $\phi\in\text{Mor}(\ca[U],\ca[W])$, and $t_0\in\bb$. Let
furthermore $H:V\to W$ be the maximal flow of $\widetilde{X}\in\ca[X](W)$, i.e.
$\pa[t]\circ H^*=H^*\circ\widetilde{X}$, subject to the initial condition
$H(t_0,\cdot)=\widetilde{\phi}:U\to W$.
Let now $\ca[V]$ be $(V,\ca[O]_{\sbb\times\ca[U]}|_{V})$ and $(t,\tau)$ the
canonical coordinates on $\sbb$, then there exists a unique $F:\ca[V]\to\ca[W]$
such that 
\begin{eqnarray}
(\text{inj }_{\bb}^{\sbb})^*\circ(\pa[t]+\pa[\tau])\circ F^* &=&  (\text{inj
}_{\bb}^{\sbb})^*\circ F^*\circ X\text{ and }\label{eqflow2}\\
F\circ\text{inj }_{\{t_0\}\times \ca[U]}^{\ca[V]} &=& \phi.\label{eqic}
\end{eqnarray}
Moreover, $\widetilde{F}:V\to W$ equals the underlying classical flow map $H$ of
the vector field $\widetilde{X}$ with initial condition $\widetilde{\phi}$. 
\end{lem}

\textbf{Proof.} Let $(u_i)=(x_i,\xi_r)$ and $(w_j)=(y_j,\eta_s)$ denote the
canonical coordinates on $\sbb[m|n]$ and $\sbb[p|q]$, respectively. Then there
exist smooth functions $a^j_J\in\ca[C]^\infty_{\bb^p}(W)$ such that
\[
X=\displaystyle\sum_{j=1}^{p+q}\left(\sum_{J}a_J^j(y)\eta^J\right)\partial_{w_j}
,
\]
where $J=(\beta_1,\ldots,\beta_q)$ runs over the index set $\{0,1\}^q$ and
$\eta^J=\displaystyle\prod_{s=1}^{q}\eta_s^{\beta_s}$.
We then have, of course,

$
X_0=\displaystyle\sum_{j}\left(\sum_{|J|=|w_j|}a_J^j(y)\eta^J\right)\partial_{
w_j}
\text{ resp. }
X_1=\displaystyle\sum_{j}\left(\sum_{|J|=|w_j|+1}a_J^j(y)\eta^J\right)\partial_{
w_j}.
$\\
Here, $|J|$ equals $\beta_1+\cdots+\beta_q\text{ mod }2$ and $|w_j|$ is the parity of the coordinate function $w_j$.
The morphism $F$ determines and is uniquely determined by functions
$f^j_I,g^j_I\in\ca[C]^\infty_{\bb\times\bb^m}(V)$ fulfilling for each
$j\in\{1,\ldots,p+q\}$
\[
F^*(w_j)=\displaystyle\sum_{|I|=|w_j|}f_I^j(t,x)\xi^I+\sum_{|I|=|w_j|+1}g_I^j(t,
x)\tau\xi^I
\]
(and $f^j_I=0$ if $|I|\neq|w_j|$, $g^j_I=0$ if $|I|\neq|w_j|+1$) as is
well-known from the standard theory of supermanifolds (compare, e.g., Thm. 4.3.1
in \cite{var}).
Here and in the sequel $I=(\alpha_1,\ldots,\alpha_n)$ is an element of the set
$\{0,1\}^n$ and $\xi^I$ stands for the product
\mbox{$\xi_1^{\alpha_1}\cdot\xi_2^{\alpha_2}\cdots \xi_n^{\alpha_n} $}. The notation $|I|$ again denotes the parity of
$I$, i.e. $|I|=\alpha_1+\cdots+\alpha_n\text{ mod }2$.\\

Equation (\ref{eqflow2}) is equivalent to the following
equations:
\begin{eqnarray}
(\text{inj }_{\bb}^{\sbb})^*\circ\pa[t]\circ F^* &=&  (\text{inj
}_{\bb}^{\sbb})^*\circ F^*\circ X_0\label{eqflow0}\\
(\text{inj }_{\bb}^{\sbb})^*\circ\pa[\tau]\circ F^* &=&  (\text{inj
}_{\bb}^{\sbb})^*\circ F^*\circ X_1\label{eqflow1}
\end{eqnarray}

Applying (\ref{eqflow0}) to the canonical coordinate functions on $\ca[W]$, we
get the following system, which is equivalent to (\ref{eqflow0}):
\begin{eqnarray}
\displaystyle\sum_{|I|=|w_j|}\pa[t]f_I^j\cdot\xi^I =
\displaystyle\sum_{|J|=|w_j|}\check{F}^*(a_J^j)\check{F}^*(\eta^J)\label{eq25}
\text{ for all } j\text{ in }\{1,\ldots,p+q\},
\end{eqnarray}
and (\ref{eqflow1}) is equivalent to 
\begin{eqnarray}
\displaystyle\sum_{|I|=|w_j|+1}g_I^j\cdot\xi^I =
\displaystyle\sum_{|J|=|w_j|+1}\check{F}^*(a_J^j)\check{F}^*(\eta^J)\label{eq26}
\text{ for all } j\text{ in }\{1,\ldots,p+q\},
\end{eqnarray}
where $\check{F}:=F\circ\text{inj
}_{\bb}^{\sbb}:\check{\ca[V]}:=(V,\ca[O]_{\bb\times\ca[U]}|_V)\to\ca[W]$. Let us
immediately observe that the underlying smooth map of $\check{F}$ equals
$\tilde{F}$, the smooth map underlying the morphism $F$.\\

Moreover the initial condition (\ref{eqic}) is equivalent to
\begin{eqnarray}
\displaystyle\sum_{|I|=|w_j|}f_I^j(t_0,x)\xi^I =\phi^*(w_j)\text{ for all } j\text{ in }\{1,\ldots,p+q\}.\label{eq27}
\end{eqnarray}
We are going to show that (\ref{eq25}) and (\ref{eq27}) uniquely determine the
functions $f^j_I$ on $V$, i.e. the morphism $\check{F}$. Then the functions
$g^j_I$ are unambiguously given by (\ref{eq26}) on $V$, and the morphism
$F$ is fully determined. \\

Let us develop Equation (\ref{eq25}) for a fixed $j$:
\begin{eqnarray}
\displaystyle\sum_{|I|=|w_j|}\pa[t]f_I^j\cdot\xi^I &=&
\displaystyle\sum_{\substack{|J|=|w_j| \\ J=(\beta_1,\dots,\beta_q)
}}\check{F}^*(a_J^j)\prod_{s=1}^{q}\check{F}^*(\eta_s)^{\beta_s}\\
\hbox{and thus }\displaystyle\sum_{|I|=|w_j|}\pa[t]f_I^j\cdot\xi^I        	
								 &=&
\displaystyle\sum_{\substack{|J|=|w_j| \\ J=(\beta_1,\dots,\beta_q)
}}\check{F}^*(a_J^j)\prod_{s=1}^{q}\left(\sum_{|L|=1}f^{p+s}_L\xi^L\right)^{
\beta_s}\label{eq29}.
\end{eqnarray}
For fixed $j$ this is an equation of Grassmann algebra-valued maps in the  
variables $t$ and $x$ that can be split in a system of scalar equations as
follows. For $K=(\alpha_1,\ldots,\alpha_n)\in\{0,1\}^n$, we will denote the
coefficient $h_K$ in front of $\xi^K$ of a superfunction $h=\sum_M
h_M(t,x)\xi^M\in\ca[O]_{\sbb[m+1|n]}$ compactly by $(h|\xi^K)$ in the sequel of
this proof.\\

Let us first describe the coefficients for $\check{F}^*(a_J^j)$ in (\ref{eq29}):
\begin{equation}
(\check{F}^*(a_J^j)|\xi^K)=0\:\:\text{ if }|K|=1\label{keyeq1}
\end{equation}
and, if $|K|=0$, 
$$(\check{F}^*(a_J^j)|\xi^K)=a_J^j\circ\widetilde{F}\:\:\text{ if
}K=(0,\ldots,0),$$
and
\begin{equation}
 (\check{F}^*(a_J^j)|\xi^K)=\displaystyle\sum_{\mu=1}^p(\partial_{y_\mu}
a^j_J)(\widetilde{F}(t,x))\cdot f^\mu_K+R\left(a^j_J,(f^\nu_I)_{\nu,
\text{deg}(I)<\text{deg}(K)}\right)\label{keyeq}
\end{equation}
$$\hspace{8cm}\text{ if deg}(K)>0.$$
Here for $I=(\alpha_1,\dots,\alpha_n)$, deg$(I)=\alpha_1+\cdots+\alpha_n$, and
-more importantly- $R=R_{j,J,K}$ is a polynomial function in $a_J^j$ and its
derivatives in the $y$-variables up to order $q$ included, and in the functions
$\{f^\nu_I|1\leq\nu\leq p+q, 0\leq \text{deg}(I)<\text{deg}(K)\}$. Equation
(\ref{keyeq1}) is obvious since $a^j_J$ is an even function, whereas equation
(\ref{keyeq}) can be deduced from standard analysis on superdomains. More
precisely, let $a$ be a smooth function on $\bb^p$ and
$\psi:\sbb[m+1|n]\to\sbb[p|q]$ morphism (of course to be applied to
$a=a^j_J,\psi=\check{F}$). Then we can develop $\psi^*(a)$ as follows (compare
the proof of Theorem 4.3.1 in \cite{var}):
\begin{eqnarray*}
\psi^*(a)&=&\displaystyle\sum_\gamma\frac{1}{\gamma!}(\partial_\gamma a )
(\widetilde{\psi}^*(y_1),\ldots,\widetilde{\psi}^*(y_p))\cdot
\prod_{\mu=1}^p(\psi^*(y_\mu)-\widetilde{\psi}^*(y_\mu))^{\gamma_\mu}\\
 &=&
a(\widetilde{\psi}(t,x))+\sum^p_{\mu=1}(\pa[y_\mu]a)(\widetilde{\psi}(t,
x))\cdot\left(\sum_{M\neq 0}f^\mu_M\cdot\xi^M\right)+
 \end{eqnarray*}
$$
\frac{1}{2}\sum_{\mu',\mu''=1}^p(\pa[y_{\mu'}]\pa[y_{\mu''}]
a)(\widetilde{\psi}(t,x))\cdot\left(\sum_{M'\neq
0}f^{\mu'}_{M'}\cdot\xi^{M'}\right)\cdot\left(\sum_{M''\neq
0}f^{\mu''}_{M''}\cdot\xi^{M''}\right)+\ldots,
$$
where $\dis\sum_{M\neq
0}f^\mu_M\cdot\xi^M=\psi^*(y_\mu)-\widetilde{\psi}^*(y_\mu)$ with $f^\mu_M$
depending on $t$ and $x$. We observe that the last RHS is a finite sum since we
work in the framework of finite-dimensional supermanifolds.\\
In order to get a contribution to $(\psi^*(a)|\xi^K)$ we can either extract
$f^\mu_K$ from the ``linear term'' or from products coming from the higher order
terms in the above development. Thus 
$$(\psi^*(a)|\xi^K)=\sum^p_{\mu=1}(\pa[y_\mu]a)(\widetilde{\psi}(t,x))\cdot
f^\mu_K+R(a,(f_I^\nu)_{\nu,\text{deg}(I)<\text{deg}(K)}),$$
where $R$ is a polynomial as described after Equation (\ref{keyeq}).\\

Furthermore, for an element $J=(\beta_1,\dots,\beta_q)$ with $|J|=0$ we have for
deg$(K)>0$
\begin{equation}
\left(\left.\prod_{s=1}^{q}\left(\sum_{|L|=1}f^{p+s}_L\xi^L\right)^{\beta_s}
\right|\xi^K\right)=R\left((f^j_I)_{j,\text{deg}(I)<\text{deg}(K)}\right)
.\label{eq312}
\end{equation}
And for an element $J=(\beta_1,\dots,\beta_q)$ with $|J|=1$ we get for
deg$(K)>0$
\begin{equation}
\left(\prod_{s=1}^{q}\left(\left.\sum_{|L|=1}f^{p+s}_L\xi^L\right)^{\beta_s}
\right|\xi^K\right)=
\left\lbrace
\begin{array}{l}
f^{p+l}_K + R\left((f^j_I)_{j,\text{deg}(I)<\text{deg}(K)}\right)\\
\text{if deg}(J)=1\text{ and }l\in\{1,\ldots,q\}\\
\text{such that }\beta_s=\delta_{s,l}\quad\forall s,\\
R\left((f^j_I)_{j,\text{deg}(I)<\text{deg}(K)}\right)\\
\text{if deg}(J)>1.
\end{array}\right.\label{eq15}
\end{equation}
Obviously, the coefficient of $\xi^K$ of the LHS of Equation (\ref{eq25}) is
given by
\[
\left(\left.\displaystyle\sum_{|I|=|w_j|}\pa[t]f_I^j\cdot\xi^I
\right|\xi^K\right)=
\left\lbrace
\begin{array}{cl}
\pa[t]f^j_K & \text{ if }|K|=|w_j|\\
0 & \text{ if }|K|=|w_j|+1
\end{array}\right.\quad\text{for } 1\leq j \leq p+q.
\]\\

Taking into account the above descriptions of the $\xi^K$-coefficients, we will
show the existence (and uniqueness) of the solution functions $\{f^j_I|1\leq
j\leq p+q,I\in\{0,1\}^n\}$ for $(t,x)\in V$ by induction on deg$(I)$ and upon
observing that all ordinary differential equations occuring are (inhomogeneous)
linear equations for the unknown functions.\\

Let us start with $\text{deg}(I)=0$ that is $I=(0,\cdots,0)$. The ``$0$-level''
of the equations (\ref{eq29}) and (\ref{eq27}) is 
$
\pa[t]f^j_{(0,\cdots,0)}=a^j_{(0,\cdots,0)}\circ\widetilde{F}\text{ and
}f_{(0,\cdots,0)}^j(t_0,x) = y_j\circ\widetilde{\phi}(x)\text{ for all }j\text{
such that }|w_j|=0.
$  
We remark that $f^j_{(0,\cdots,0)}$ is simply $y_j\circ\widetilde{F}$ and
$a^j_{(0,\cdots,0)}$ is $\widetilde{X}(y_j)$. Thus $\widetilde{F}$ is the flow
of $\widetilde{X}$ with initial condition $\widetilde{\phi}$ at $t=t_0$, i.e.,
$\widetilde{F}=H$ on $V$.
Thus the claim is true for $I=(0,\cdots,0)$.\\

 Suppose $k>0$ and that the functions $f^j_I$ are uniquely defined on $V$ for
all $j$ and all $I$ such that $\text{deg}(I)< k$.
Let $K$ be such that $\text{deg}(K)=k$. Let us distinguish the two possible
parities of $k$ in order to determine $f^j_K$ for all $j$. Recall that $f^j_K=0$
if the parities of $K$ and $j$ are different.\\

If $k$ is even, i.e., $|K|=0$, we only have to consider $j$ such that $|w_j|=0$.
Putting (\ref{keyeq}) and (\ref{eq312}) together, we find in this case\\

$\begin{array}{ccl}
\pa[t]f_K^j &=&\left(\left.\displaystyle\sum_{\substack{|J|=0 \\
J=(\beta_1,\dots,\beta_q)
}}\check{F}^*(a_J^j)\prod_{s=1}^{q}\left(\sum_{|L|=1}f^{p+s}_L\xi^L\right)^{
\beta_s}\right|\xi^K\right)\vspace*{2mm}\\
&=&\left(\left.\displaystyle\sum_{\substack{\text{deg}(J)=0 \\
J=(\beta_1,\dots,\beta_q)
}}\check{F}^*(a_J^j)\prod_{s=1}^{q}\left(\sum_{|L|=1}f^{p+s}_L\xi^L\right)^{
\beta_s}\right.\right.
\end{array}$\vspace*{2mm}\\
$\hspace*{3cm}\left.\left.+\displaystyle\sum_{\substack{|J|=0 \\ \text{deg}(J)>0 \\
J=(\beta_1,\dots,\beta_q)
}}\check{F}^*(a_J^j)\prod_{s=1}^{q}\left(\sum_{|L|=1}f^{p+s}_L\xi^L\right)^{
\beta_s}\right|\xi^K\right)$\vspace{4mm}\\
$=\left(\left.\displaystyle\check{F}^*(a_{(0,\cdots,0)}^j)+\displaystyle\sum_{
\substack{|J|=0 \\ \text{deg}(J)>0 \\ J=(\beta_1,\dots,\beta_q)
}}\check{F}^*(a_J^j)\prod_{s=1}^{q}\left(\sum_{|L|=1}f^{p+s}_L\xi^L\right)^{
\beta_s}\right|\xi^K\right)$\vspace*{2mm}\\
$=\displaystyle\sum_{\mu=1}^p\left(\partial_{y_\mu}a^j_{(0,\cdots,0)}\circ
\widetilde{F}\right) f^\mu_K+R\left((a^j_J)_{J},(f^\nu_I)_{\nu,
\text{deg}(I)<\text{deg}(K)}\right).$\\

Moreover, the initial condition gives
$
\displaystyle f_K^j(t_0,x) = (\phi^*(y_j)|\xi^K),$ for all $j\text{ in }\{1,\ldots,p\}.
$
Since the $a^j_J$ are the (given) coefficients of the vector field $X$ and the
functions $f^\nu_I$ with deg$(I)<k$ are known by the induction
hypothesis, we have a unique local solution function $f^j_K$. Since the ordinary
differential equation for $f^j_K$ is linear its solution is already defined for
all $(t,x)\in V$. Thus in the case that $k$ is even $f^j_K$ is unambiguously
defined on $V$ for all $j\in\{1,\ldots,p+q\}$ and for all $K$ with
deg$(K)=k$.\\ 

Now, if $k$ is odd, i.e., $|K|=1$, we only have to consider $j$ such that
$|w_j|=1$. Using (\ref{keyeq}) and (\ref{eq15}), we find in this case:\\

$\begin{array}{ccl}
\pa[t]f_K^j &=&\left(\left.\displaystyle\sum_{\substack{|J|=1 \\
J=(\beta_1,\dots,\beta_q)
}}\check{F}^*(a_J^j)\prod_{s=1}^{q}\left(\sum_{|L|=1}f^{p+s}_L\xi^L\right)^{
\beta_s}\right|\xi^K\right)\vspace*{2mm}\\
&=&\left(\left.\displaystyle\sum_{\substack{\text{deg}(J)=1 \\
J=(\beta_1,\dots,\beta_q)
}}\check{F}^*(a_J^j)\prod_{s=1}^{q}\left(\sum_{|L|=1}f^{p+s}_L\xi^L\right)^{
\beta_s}+\right.\right.
 \end{array}$\vspace*{2mm}\\
$\left.\left.\displaystyle\sum_{\substack{|J|=1 \\ \text{deg}(J)>1 \\
J=(\beta_1,\dots,\beta_q)
}}\check{F}^*(a_J^j)\prod_{s=1}^{q}\left(\sum_{|L|=1}f^{p+s}_L\xi^L\right)^{
\beta_s}\right|\xi^K\right)$\vspace{4mm}\\
$=\left(\left.\displaystyle\sum_{s=1}^q\check{F}^*(a_{(\delta_{1s},\cdots,\delta_{
qs})}^j)\left(\sum_{|L|=1}f^{p+s}_L\xi^L\right)\right.\right.$

$\hspace*{2cm}\left.\left.+\displaystyle\sum_{\substack{
|J|=1 \\ \text{deg}(J)>1\\ J=(\beta_1,\dots,\beta_q)
}}\check{F}^*(a_J^j)\prod_{s=1}^{q}\left(\sum_{|L|=1}f^{p+s}_L\xi^L\right)^{
\beta_s}\right|\xi^K\right)$\vspace{4mm}\\
$= \displaystyle\sum_{s=1}^q
\left(a^j_{(\delta_{1s},\cdots,\delta_{qs})}\circ
\widetilde{F}\right)f^{p+s}_K+R\left((a^j_J)_{J},(f^\nu_I)_{\nu,
\text{deg}(I)<\text{deg}(K)}\right).$\\

Moreover, the initial condition gives
\begin{eqnarray*}
\displaystyle f_K^j(t_0,x) = (\phi^*(w_j)|\xi^K)\text{ for all }
j\text{ in }\{p+1,\ldots,p+q\}.
\end{eqnarray*}
It follows as in the case of $|K|=0$, that $f^j_K$ exists uniquely for all
$(t,x)\in V$, for all $j\in\{1,\ldots,p+q\}$ and for all $K$ with
$\text{deg}(K)= k$. \\

We conclude that the functions $\{f^j_I|1\leq j\leq p+q,I\in\{0,1\}^n\}$ are
uniquely defined on the whole of $V$. Since the $\{g^j_I|1\leq j\leq
p+q,I\in\{0,1\}^n\}$ are determined by Equation (\ref{eq26}) from the
$\{f^j_I|1\leq j\leq p+q\}$ via comparison of coefficients, the morphism
$F:\ca[V]\to\ca[W]$ is uniquely determined.
\hfill$\square$\\

We now consider the global problem of integrating a vector field on a
supermanifold. In order to prove that there exists a unique maximal flow of a
vector field, the following lemma will be crucial.
\begin{lem}\label{theoglob}Let $\cat$ and $\cat[S]$ be supermanifolds, $X$ a
vector field in $\tm(M)$ and $\phi$ in $\text{Mor}(\ca[S],\ca)$. Then
\begin{enumerate}[(i)]
 \item there exists an open sub supermanifold $\ca[V]=(V,\ca[O]_{\sbb\times \ca[S]}|_{V})$ of $\sbb\times\ca[S]$ with
$V$
open in $\bb\times S$ such  that $\{0\}\times S\subset V$ and for all $ x$ in $ S$, $ (\bb\times\{x\})\cap V$ is an
interval, and a morphism $F:\ca[V]\to\ca$
satisfying:
\begin{eqnarray}
(\text{inj }_{\bb}^{\sbb})^*\circ(\pa[t]+\pa[\tau])\circ F^* &=&  (\text{inj
}_{\bb}^{\sbb})^*\circ F^*\circ X\text{ and }\label{eq2110}\\
F\circ\text{inj }_{\{0\}\times \ca[S]}^{\ca[V]} &=& \phi\:.\label{eq211}
\end{eqnarray}
 \item Let furthermore $F_1:\ca[V]_1\to\ca$ and $F_2:\ca[V]_2\to\ca$ be morphisms satisfying (\ref{eq2110}) and
(\ref{eq211})
where $\ca[V]_i=(V_i,\ca[O]_{\sbb\times\ca[S]}|_{V_i})$ with $V_i$ open in $\bb\times
S$ such that $\{0\}\times S\subset V_i$ and for all $x$ in $S$,
$(\bb\times\{x\})\cap V_i$ is an interval, for $i=1,2$.
Then $F_{1|\ca[V]_{12}}=F_{2|\ca[V]_{12}}$ on
$\ca[V]_{12}=(V_{12},\ca[O]_{\sbb\times\ca[S]}|_{V_{12}})$, where $V_{12}=V_{1}\cap V_{2}$.
\end{enumerate}
\end{lem}
\textbf{Proof.} (i) Let $\widetilde{\phi}:S\to M$ denote the induced map of the underlying
classical manifolds. Given now $s$ in $S$ and coordinate domains $\ca[U]_s$ of $s$ and $\ca[W]_s$ of
$\widetilde{\phi}(s)$, isomorphic to
superdomains $\check{\ca[U]}_s\subset\sbbmn$ resp. $\check{\ca[W]}_s\subset\sbbpq$, by Lemma \ref{theolocal} we
get solutions of (\ref{eq2110}) and (\ref{eq211}) near $s$ (upon reducing the size of $\ca[U]_s$ if necessary):
$\sbb\times\ca[S]\supset\sbb\times \ca[U]_s\supset\ca[V]_s \overset{F^s}{\longrightarrow}\ca[W]_s\subset\ca. $
If $\ca[V]_{s_1}\cap \ca[V]_{s_2}\neq \emptyset$ (compare Figure 1) we know, again by Lemma
\ref{theolocal}, that $F^{s_1}$ and $F^{s_2}$
coincide on this intersection. Thus, by taking the union $\ca[V]$ of $\ca[V]_{s}$ for all $s$ in $S$, we get a morphism
$F:\sbb\times \ca[S]\supset\ca[V]\to \ca$ such that $F_{|\ca[V]_s}=F^s$ for all $s$, and fulfilling (\ref{eq2110}) and
(\ref{eq211}).
\begin{center}
  \includegraphics[scale=.6]{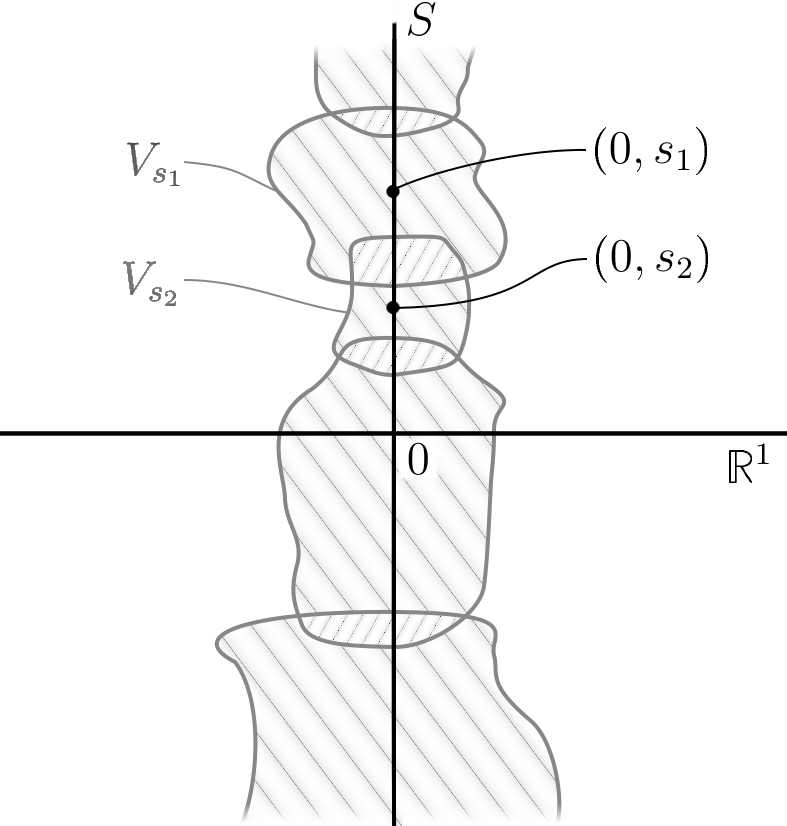}\\
\textsc{Figure 1}
\end{center}
(ii) We define $A$ as the set of points $(t,x)\in
V_{12}$ such that
there exists $\epsilon=\epsilon_{(t,x)}>0$ and $\ca[U]=\ca[U]
_{(t,x)}$ an open sub supermanifold of $\ca[S]$, such that its body $U$ contains
$x$ and for $\ca[V]=\ca[V]_{(t,x)}=(V_{(t,x)},\ca[O]_{\sbb\times\ca[S]})=(]-\epsilon,t+\epsilon[\times
U,\ca[O]_{\sbb\times\ca[S]})$ we have $F_{1|\ca[V]}=F_{2|\ca[V]}$. Of course, if $t<0$ the interval will be of the type
$]t-\epsilon, \epsilon[$ (See Figure 2). The claim of
the Lemma is now equivalent to $A=V_{12}$. The set $A$ is obviously open.
\begin{center}
\includegraphics[scale=.53]{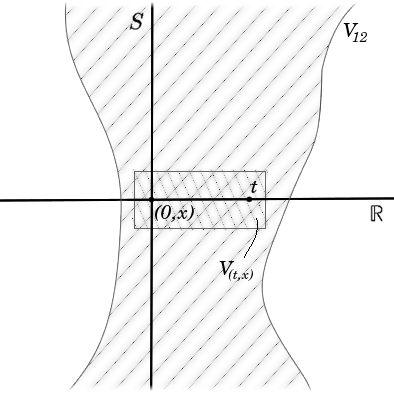}\\
\textsc{Figure 2}
\end{center}
By an easy application of Lemma \ref{theolocal}, A contains $\{0\}\times S$.
The assumptions imply that for all $ x\in S$, the set $I_x\subset \bb$, defined by $
(\bb\times\{x\})\cap V_{12}=I_x\times\{x\}$, is an open interval containing $0$. The definition of $A$ implies that the
set
$J_x\subset I_x$ such that $(\bb\times\{x\})\cap A=J_x\times\{x\}$ is an open interval containing $0$ as well.
Assuming now that $A\neq V_{12}$, then there exists a point $(t,x_0)\in
V_{12}\backslash A$ such that $J_{x_0}\neq I_{x_0}$. Without loss of generality we can assume that $t>0$ and that for
$0\leq t'<t$, $(t',x_0)\in A$.
Let $U_0$ be an
open coordinate neighborhood of $x_0$ in $S$ and $\delta>0$ such that, with
$V_0:=]t-\delta,t+\delta[\times U_0\subset V_{12}$,
 $H(V_0)\subset W$, where $\ca[W]=(W,\om|_{W})$ is a coordinate patch of $\ca$.
Choose $t_0\in ]t-\delta,t[$. Then $(t_0,x_0)\in A$ and thus there exists
$\epsilon>0$ and $\ca[U]$ an open sub supermanifold of $\ca[U]_0=(U_0,\om[S]|_{U_0})$
containing $x_0$ such that 
\begin{equation}
F_{1|\ca[V]}=F_{2|\ca[V]},\:\text{ where }
\ca[V]=]-\epsilon,t_0+\epsilon[\times\sbb[0|1]\times\ca[U]\subset\ca[V]_{12}
.\label{eqlemma32}
\end{equation}
On $\ca[V]'=]t-\delta,t+\delta[\times\sbb[0|1]\times\ca[U]\subset \ca[V]_{12}$, $F_1$ and $F_2$ are
defined and for
$i=1,2$ the maps $F_i\circ\text{inj}^{\ca[V]}_{\{t_0\}\times\ca[U]}$ coincide by (\ref{eqlemma32}) (Compare Figure
3 for the relative positions of the underlying topological spaces of these open sub supermanifolds of
$\sbb\times\ca[S]$).
\begin{center}
\includegraphics[scale=.95]{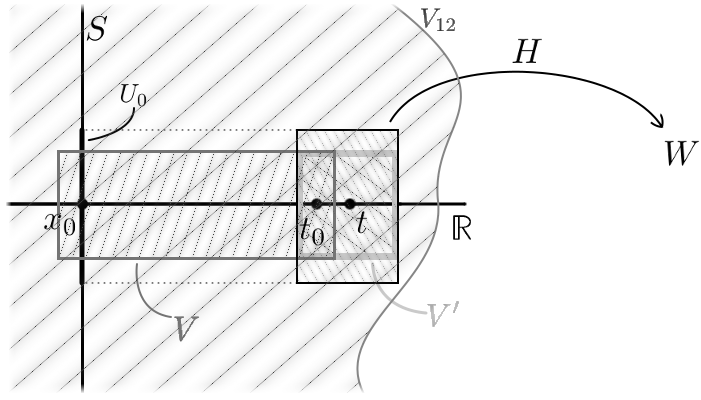}\\
\textsc{Figure 3}
\end{center}
By Lemma
\ref{theolocal} we have $F_{1|\ca[V]'}=F_{2|\ca[V]'}$. Thus $F_1=F_2$ on $\ca[V]\cup\ca[V]'$, and we conclude that
$(t,x_0)\in A$. This contradiction shows that $V_{12}=A$.\hfill$\square$\\

\textbf{Remarks. }(1) Obviously, Lemma \ref{theoglob} holds true for an arbitrary $t_0\in\bb$ replacing $t_0=0$.\\
(2) Let us call a ``flow domain for $X$ with initial condition $\phi\in\text{Mor}(\ca[S],\ca)$ (with respect to
$t_0\in\bb$)'' a domain $\ca[V]\subset\sbb\times\ca[S]$ such that $\{t_0\}\times S\subset V$ and for all $s $ in $S$,
$(\bb\times\{s\})\cap V$ is connected, i.e. an interval (times $\{s\}$) and such that a solution $F$ (a ``flow'') of
(\ref{eq2110}) and (\ref{eq211}) exists on $\ca[V]$. By the preceding lemma there exists such ``flow domains''.

\begin{theo}\label{theoflowm}
Let $\ca$ and $\ca[S]$ be supermanifolds, $X$ be a vector field in
$\tm(M)$, $\phi$ in $\text{Mor}(\ca[S],\ca)$ and $t_0$ in $\bb$. Then there exists a
unique map $F:\ca[V]\to\ca$ such that
\begin{eqnarray*}
(\text{inj }_{\bb}^{\sbb})^*\circ(\pa[t]+\pa[\tau])\circ F^* &=&  (\text{inj
}_{\bb}^{\sbb})^*\circ F^*\circ X\text{ and }\\
F\circ\text{inj }_{\{t_0\}\times \ca[S]}^{\ca[V]} &=& \phi\:,
\end{eqnarray*}
where $\ca[V]=(V,\ca[O]_{\sbb\times\ca[U]}|_V)$ is the maximal flow domain for $X$
with the given initial condition.\\
Moreover, $\widetilde{F}:V\to M$ is the maximal flow of
$\widetilde{X}\in\ca[X](M)$ subject to the initial condition $\widetilde{\phi}$
at $t=t_0$.
\end{theo}
\textbf{Proof.} The proof of the theorem follows immediately from the Lemmata
\ref{theolocal} and \ref{theoglob}, upon taking the union of all flow domains and flows for $X$ as defined in the
preceding remark.\hfill$\square$

\section{Supervector fields and local $\sbb$-actions}
Given a vector field on a classical, ungraded, manifold, the flow map
$\widetilde{F}$ (for $S=M$, $\widetilde{\phi}=\text{id}_{M}$) is always a local
action of $\bb$ with its usual (and unique up to isomorphism) Lie group structure, the standard addition. The flow maps
for vector fields
described in the preceding section (taking here $\ca[S]=\ca$,
$\phi=\text{id}_{\ca}$), do not always have the analogous property of being
local actions of $\sbb$ with an appropriate Lie supergroup structure. Two
characterizations of those vector fields $X=X_{0}+X_{1}$ that generate a
local $\sbb$-action were found by J. Monterde and O. A. S\'anchez-Valenzuela. We
will give in this section a short proof of a slightly more general result, whose condition (iii)
seems to be more easily verified in practice than those given in \cite{Mo-SV2}
(compare Thm. 3.6 and its proof there).\\

Let us begin by giving a useful two-parameter family of Lie supergroup structures on the supermanifold $\sbb$ and their
right invariant vector fields.\\

\begin{lem}\label{rivf}
Let $a$ and $b$ be real numbers such that $a\cdot b=0$ and
$\mu_{a,b}=\mu:\sbb\times\sbb\to\sbb$ be defined by
\begin{eqnarray*}
\widetilde{\mu}(t_1,t_2)&=&t_1+t_2,\\
\mu^*(t)&=& t_1+t_2+a\tau_1\tau_2, \\
\mu^*(\tau)&=& \tau_1+e^{bt_1}\tau_2\qquad.
\end{eqnarray*}
Then \begin{enumerate}[(i)]
\item there exists a unique Lie supergroup structure on $\sbb$ such that the
multiplication morphism is given by $\mu_{a,b}$,
\item the right invariant vector fields on $(\sbb,\mu_{a,b})$ are given by
the graded vector space $\bb D_0\oplus \bb D_1$, where
$$D_0:=\pa[t]+b\cdot\tau\pa[\tau]\text{ and }D_1:=\pa[\tau]+a\cdot\tau\pa[t],$$
and they obey $[D_0,D_0]=0$, $[D_0,D_1]=-bD_1$ and $[D_1,D_1]=2aD_0$.
\end{enumerate}
\end{lem}
\textbf{Proof. }  The first assertion follows by an easy direct verification.\\
 Let now $Z$ be a vector field on $\sbb$, then there exist
$\alpha,\beta,\gamma,\eta\in\ca[C]^\infty_{\bb}(\bb)$ such that $Z=Z_0+Z_1$
where $Z_0=\alpha\pa[t]+\beta\tau\pa[\tau]$ and
$Z_1=\gamma\pa[\tau]+\eta\tau\pa[t]$.\\
  Recall that $Z$ is a right invariant vector field if $Z$ fulfills the
following condition:
$$
\mu^*\circ Z=(Z\otimes \text{id})\circ \mu^*.
$$
Remark that this condition can be split into the following two conditions:
$$
\mu^*\circ Z_0=(Z_0\otimes \text{id})\circ \mu^*\text{ and }\mu^*\circ
Z_1=(Z_1\otimes \text{id})\circ \mu^*.
$$
Thus, if $Z$ is a right invariant vector field, the functions
$\alpha,\beta,\gamma,\eta$ are characterized as follows:
\begin{eqnarray*}
\mu^*\circ
Z_0(t)&=&\mu^*(\alpha)=\alpha(t_1+t_2)+a\alpha'(t_1+t_2)\tau_1\tau_2,\\
\\
(Z_0\otimes\text{id})\circ\mu^*(t)&=&(Z_0\otimes\text{id})(t_1\otimes 1+1\otimes
t_2+ a\tau_1\otimes\tau_2)\\
&=& Z_0(t_1)\otimes 1+Z_0(1)\otimes t_2+ aZ_0(\tau_1)\otimes\tau_2\\
&=& \alpha(t_1)+a\beta(t_1)\tau_1\tau_2\qquad.\\
\end{eqnarray*}
Consequently $\alpha$ is a constant and $a\cdot\beta=0$. Similarly,
\begin{eqnarray*}
\mu^*\circ
Z_0(\tau)&=&\mu^*(\beta\tau)\\
&=&\left(\beta(t_1+t_2)+a\beta'(t_1+t_2)\tau_1
\tau_2\right)\cdot(\tau_1+e^{bt_1}\tau_2)\\
&=&\beta(t_1+t_2)\tau_1+\beta(t_1+t_2)e^{bt_1}\tau_2,\\
\\
(Z_0\otimes\text{id})\circ\mu^*(\tau)&=&(Z_0\otimes\text{id})(\tau_1\otimes
1+e^{bt_1}\otimes\tau_2)\\
&=& Z_0(\tau_1)\otimes 1+Z_0(e^{bt_1})\otimes \tau_2\\
&=& \beta(t_1)\tau_1+\alpha\cdot be^{bt_1}\tau_2,
\end{eqnarray*}
implying that $\beta$ is a constant and $b\cdot\alpha=\beta$. For the odd part
we find
\begin{eqnarray*}
\mu^*\circ
Z_1(t)&=&\mu^*(\eta\tau)=\eta(t_1+t_2)\tau_1+\eta(t_1+t_2)e^{bt_1}\tau_2,\\
\\
(Z_1\otimes\text{id})\circ\mu^*(t)&=&(Z_1\otimes\text{id})(t_1\otimes 1+1\otimes
t_2+ a\tau_1\otimes\tau_2)\\
&=&\eta(t_1)\tau_1+a\gamma(t_1)\tau_2.\\
\end{eqnarray*}
Thus $\eta$ is a constant and $a\cdot e^{-bt_1}\cdot\gamma=\eta$. Finally,
\begin{eqnarray*}
\mu^*\circ
Z_1(\tau)&=&\mu^*(\gamma)=\gamma(t_1+t_2)+a\gamma'(t_1+t_2)\tau_1\tau_2,\\
\\
(Z_1\otimes\text{id})\circ\mu^*(\tau)&=&\gamma(t_1)+\eta(t_1)be^{bt_1}
\tau_1\tau_2,\\
\end{eqnarray*}
thus $\gamma$ is a constant and $b\cdot\eta=0$.\\

 Since $a\cdot b=0$, it suffices to consider the cases $a=0$ and $b=0$. If
$a=0$, we have $Z_0=\alpha\pa[t]+b\alpha\tau\pa[\tau]=\alpha D_0$ and
$Z_1=\gamma\pa[t]=\gamma D_1$ where $\alpha,\gamma\in\bb$. If $b=0$, we have
$Z_0=\alpha\pa[t]=\alpha D_0$ and $Z_1=\gamma\pa[\tau]+a\gamma\tau\pa[t]=\gamma
D_1$ where $\alpha,\gamma$ are real numbers. Reciprocally, in both cases, $Z$ is
a right invariant vector field.
\hfill$\square$\\

\textbf{Remarks. }(1) It can easily be checked that the above family yields only three non-isomorphic Lie supergroup
structures on $\sbb$, since $(\sbb, \mu_{a,0})$ with $a\neq0$ is isomorphic to $(\sbb,\mu_{1,0})$ and $(\sbb,\mu_{0,b})$
with $b\neq0$ is isomorphic to $(\sbb,\mu_{0,1})$ and the three multiplications $\mu_{0,0}$, $\mu_{1,0}$ and
$\mu_{0,1}$ correspond to non-isomorphic Lie supergroup structures on $\sbb$. Nevertheless it is very convenient to
work here with the more flexible two-parameter family of multiplications.\\
(2) In fact, all Lie supergroup structures on $\sbb$ are equivalent to $\mu_{0,0}$, $\mu_{1,0}$ or $\mu_{0,1}$. See,
e.g., \cite{F-O} for a direct approach to the classification of all Lie supergroup structures on $\sbb$.

\begin{de}
 Let $\ca[G]=(G,\om[G])$ resp. $\ca=(M,\om)$ be a Lie supergroup with multiplication morphism $\mu$ and unit element $e$
resp. a
 supermanifold. A ``local action of $\ca[G]$ on $\ca$'' is given by the following data:\\
 a collection $\Pi$ of pairs
of
open subsets $\pi=(U_{\pi},W_{\pi})$ of $M$, where $U_{\pi}$ is relatively compact in $W_{\pi}$, with
associated open
sub supermanifolds $\ca[U]_{\pi}\subset \ca[W]_{\pi}\subset\ca$ such that $\{U_{\pi}|\pi\in\Pi\}$ is an open
covering of $M$, and for all $\pi$ in $\Pi$ an open sub supermanifold $\ca[G]_{\pi}\subset\ca[G]$, containing the
neutral
element $e$ and a morphism
$$\Phi_{\pi}:\ca[G]_{\pi}\times\ca[U]_{\pi}\to\ca[W]_{\pi}$$ 
fulfilling 
\begin{enumerate}[(1)]
\item $\Phi_{\pi}\circ\left( e\times\text{id}_{\,\ca[U]_{\pi}} \right)=\text{id}_{\,\ca[U]_{\pi}}$, where
$e:\{\text{pt}\}\to\ca[G]$ is viewed as a morphism,
 \item $\Phi_{\pi}\circ\left( \mu\times\text{id}_{\ca} \right)=\Phi_{\pi}\circ\left(
\text{id}_{\,\ca[G]}\times\Phi_{\pi}
\right),$ where both sides are defined,
\item if $U_{\pi}\cap U_{\pi'}\neq\emptyset$, $\Phi_{\pi}=\Phi_{\pi'}$ on $\left( \ca[G]_{\pi}\cap\ca[G]_{\pi'}
\right)\times\left( \ca[U]_{\pi}\cap\ca[U]_{\pi'}
\right).$
\end{enumerate}
\end{de}

\begin{prop}\label{LocalAction}
 Let $\ca[G]=\left( G,\om[G] \right)$ resp. $\ca=\left( M,\om \right)$ be a Lie supergroup resp. a supermanifold.
Then 
\begin{enumerate}[(i)]
 \item a local $\ca[G]$-action on $\ca$, specified by a set $\Pi$ and morphisms
$\{(\ca[U]_{\pi},\ca[W]_{\pi},\ca[G]_{\pi},\Phi_{\pi})|\pi\in\Pi\},$ gives rise to an open sub supermanifold
$\ca[V]\subset\ca[G]\times\ca$ containing $\{e\}\times\ca$ and a morphism $\Phi_{\ca[V]}:\ca[V]\to\ca$ such that 
\begin{equation}
 \Phi_{\ca[V]}\circ\left( \mu\times\text{id}_{\ca} \right) = \Phi_{\ca[V]}\circ\left(
\text{id}_{\,\ca[G]}\times\Phi_{\ca[V]} \right),\tag{$\ast$}\label{actionformula}
\end{equation}
 where both sides are defined and such that
\begin{equation}
 \Phi_{\pi}=\Phi_{\ca[V]}\text{ on }\left(\ca[G]_{\pi}\times \ca[U]_{\pi}
\right)\cap\ca[V],\:\forall\pi\in\Pi,\tag{$\ast\ast$}\label{compatibformula}
\end{equation}
\item an open sub supermanifold $\ca[V]\subset\ca[G]\times\ca$ containing $\{e\}\times\ca$ and a morphism
$\Phi:\ca[V]\to\ca$
such that (\ref{actionformula}) is fulfilled, where it makes sense, yields a local $\ca[G]$-action on $\ca$ such that
(\ref{compatibformula}) holds.
\end{enumerate}
\end{prop}
\textbf{Proof.} As in the classical case of ungraded manifolds and Lie groups.\hfill$\square$

\begin{theo}\label{theoaction}
Let $\ca$ be a supermanifold, $X$ a vector field on $\ca$ and
$\ca[V]\subset\sbb\times\ca$ the domain of the maximal flow $F:\ca[V]\to\ca$
satisfying
\begin{eqnarray*}
(\text{inj }_{\bb}^{\sbb})^*\circ(\pa[t]+\pa[\tau])\circ F^* &=&  (\text{inj
}_{\bb}^{\sbb})^*\circ F^*\circ X\text{ and }\\
F\circ\text{inj }_{\{0\}\times \ca}^{\ca[V]} &=& \text{id}_{\ca}.
\end{eqnarray*}
 Let $a$ and $b$ be real numbers. Then the following assertions are equivalent:
\begin{enumerate}[(i)]
\item the map $F$ fulfills
\[
(\pa[t]+\pa[\tau]+\tau(a\pa[t]+b\pa[\tau]))\circ F^*=F^*\circ X,
\]
\item the map $F$ is a local $(\sbb,\mu_{a,b})$-action on $\ca$,
\item $\bb X_0\oplus\bb X_1$ is a sub Lie superalgebra of $\tm(M)$ with commutators $[X_0,X_1]=-bX_1$ and
$[X_1,X_1]=2aX_0$.
\end{enumerate}
\end{theo}

\textbf{Proof.} Recall that $F$ fulfills
\begin{eqnarray*}
(\text{inj }_{\bb}^{\sbb})^*\circ\pa[t]\circ F^* &=&  (\text{inj
}_{\bb}^{\sbb})^*\circ F^*\circ X_0\:\:\text{ and }\\
(\text{inj }_{\bb}^{\sbb})^*\circ\pa[\tau]\circ F^* &=&  (\text{inj
}_{\bb}^{\sbb})^*\circ F^*\circ X_1.
\end{eqnarray*}
 Denoting the projection from $\sbb$ to $\bb$ by $p$, we have 
 \begin{eqnarray*}
\text{id}_{\sbb}^*&=& p^*\circ(\text{inj }_{\bb}^{\sbb})^*+\tau\cdot
p^*\circ(\text{inj }_{\bb}^{\sbb})^*\circ\pa[\tau]
\end{eqnarray*}
which we will write more succinctly as
 \begin{eqnarray}
\text{id}_{\sbb}^*&=& (\text{inj }_{\bb}^{\sbb})^*+\tau\cdot (\text{inj
}_{\bb}^{\sbb})^*\circ\pa[\tau].\label{astus}
\end{eqnarray}
Using relation $(\ref{astus})$ and the equations fulfilled by $F^*$ we get
\begin{eqnarray*}
F^*\circ X &=& (\text{inj }_{\bb}^{\sbb})^*\circ F^*\circ X + \tau\cdot
(\text{inj }_{\bb}^{\sbb})^*\circ\pa[\tau]\circ F^*\circ X\\
		   &=&  (\text{inj }_{\bb}^{\sbb})^*\circ(\pa[t]+\pa[\tau])\circ
F^* + \tau\cdot (\text{inj }_{\bb}^{\sbb})^* \circ F^*\circ X_1\circ X\\
		   &=&  (\text{inj }_{\bb}^{\sbb})^*\circ(\pa[t]+\pa[\tau])\circ
F^* \\
                   &+& \tau\cdot (\text{inj }_{\bb}^{\sbb})^* \circ F^*\left([X_1,X_0]+X_0\circ
X_1+\frac{1}{2}[X_1,X_1]\right).
\end{eqnarray*}
Since
\begin{eqnarray*}
(\text{inj }_{\bb}^{\sbb})^* \circ F^*\circ X_0\circ X_1 &=& (\text{inj
}_{\bb}^{\sbb})^*\circ\pa[t] \circ F^*\circ X_1\\
&=&\pa[t]  \circ(\text{inj }_{\bb}^{\sbb})^*\circ F^*\circ X_1 \\
&=&\pa[t]  \circ(\text{inj }_{\bb}^{\sbb})^*\circ \pa[\tau]\circ F^* \\
&=&\pa[t] \circ \pa[\tau]\circ F^*\\
&=& \pa[\tau] \circ \pa[t]\circ F^*,
\end{eqnarray*}
we arrive at
\begin{eqnarray*}
F^*\circ X	   &=& (\text{inj }_{\bb}^{\sbb})^*\circ(\pa[t]+\pa[\tau])\circ
F^* +\hspace*{40mm}
\end{eqnarray*}
\vspace{-10mm}
\begin{eqnarray}
	\hspace*{2cm}	&& \tau\cdot 
F^*\left([X_1,X_0]+\frac{1}{2}[X_1,X_1]\right) + \tau\cdot
\pa[\tau]\circ\pa[t]\circ F^*.\label{astus2}
\end{eqnarray}		   
On the other hand, if $a$ and $b$ are real numbers, we have, again using
(\ref{astus})
\begin{eqnarray*}
&&(\pa[t]+\pa[\tau]+\tau(a\pa[t]+b\pa[\tau]))\circ F^* \\
&&= \left((\text{inj
}_{\bb}^{\sbb})^*+\tau\cdot (\text{inj
}_{\bb}^{\sbb})^*\circ\pa[\tau]\right)\circ(\pa[t]+\pa[\tau])\circ F^*\\
&&\hspace*{1cm} +\tau\cdot\left(a\cdot (\text{inj }_{\bb}^{\sbb})^*\circ\pa[t]+b
\cdot(\text{inj }_{\bb}^{\sbb})^*\circ\pa[\tau]\right)\circ F^*\\
&&= (\text{inj }_{\bb}^{\sbb})^*\circ(\pa[t]+\pa[\tau])\circ F^*+\tau\cdot
\pa[\tau]\circ\pa[t]\circ F^*\\
&& \hspace*{1cm}+\tau\cdot F^*\circ\left(a X_0+b X_1\right).
\end{eqnarray*}
Thus we have
\begin{equation}
 \begin{array}{c}
\dis(\pa[t]+\pa[\tau]+\tau(a\pa[t]+b\pa[\tau]))\circ F^*-F^*\circ X\hspace*{3cm}\vspace*{1mm}\\
 \hspace*{1cm}=\dis\tau\cdot
F^*\circ\left(aX_0-\frac{1}{2}[X_1,X_1]+bX_1-[X_1,X_0]\right).
\end{array}\label{eq318}
\end{equation}
Since $F$ satisfies the initial condition 
$(\text{inj }_{\{0\}\times \ca[M]}^{\ca[V]})^*\circ F^* = \text{id}_{\ca},$
$\tau\cdot F^*$ is injective and thus Equation (\ref{eq318}) easily implies the equivalence of
$(i)$ and $(iii)$.\\

We remark that, in this case, we automatically have $a\cdot b=0$ since the
Jacobi identity implies that $
 [ X_1,[ X_1, X_1]]=[[ X_1, X_1], X_1]+(-1)^{1\cdot 1}[ X_1,[ X_1, X_1]]\text{,
i.e.,  }2a\cdot b\cdot X_1=[ X_1,[ X_1, X_1]]=0$.\\

Assume now that $a$ and $b$ are real numbers such that $(i)$ satisfied, and let
$\mu=\mu_{a,b}$ be as in Lemma \ref{rivf}. We have to show that $F$ is a local
action of $(\sbb,\mu)$. \

Let us define $$G:=F\circ(\text{id}_{\sbb}\times
F):\sbb\times(\sbb\times\ca)\to\ca$$ and $$H:=F\circ(\mu\times\text{id}_{\ca
}):(\sbb\times\sbb)\times\ca\cong\sbb\times(\sbb\times\ca)\to\ca.$$
In order to prove that $F$ is a $\sbb$-action on, we have to show that $G=H$. We
observe that $G$ is the integral curve of $X$ subject to the initial condition
$F\in\text{Mor}(\sbb\times\ca,\ca)$.\\

Let us prove that the morphism $H$ satisfies the following conditions:
\begin{eqnarray}
\hspace*{-6mm}(\text{inj
}_{\bb\times(\sbb\times\ca)}^{\sbb\times(\sbb\times\ca)})^*\circ(\pa[t_1]+\pa[
\tau_1])\circ H^* &=& (\text{inj
}_{\bb\times(\sbb\times\ca)}^{\sbb\times(\sbb\times\ca)})^*\circ H^*\circ X 
\label{eq212}\\
 H\circ\text{inj }_{\{0\}\times(\sbb\times\ca)}^{\sbb\times(\sbb\times\ca)} &=&
F.\label{eq213}
\end{eqnarray}
Then by the unicity of integral curves we have $H=G$. \\

Equation (\ref{eq213}) holds true since $\mu\circ\text{inj
}_{\{0\}\times\sbb}^{\sbb\times\sbb}=\text{id}_{\sbb}$. \\

Defining $D:=D_0+D_1=\pa[t_1]+\pa[\tau_1]+\tau_1(a\pa[t_1]+b\pa[\tau_1])$ and
writing $\text{inj }_{|t_1}$ for $\text{inj
}_{\bb\times(\sbb\times\ca)}^{\sbb\times(\sbb\times\ca)}$ and using
right invariance of $D$, we arrive 
at equation ($\ref{eq212}$) as follows
\[\begin{array}{l}
\hspace*{-3cm}(\text{inj
}_{\bb\times(\sbb\times\ca)}^{\sbb\times(\sbb\times\ca)})^*\circ(\pa[t_1]+\pa[
\tau_1])\circ H^* \\
\end{array}\]
\begin{eqnarray*}
\hspace*{3cm}&=&(\text{inj }_{|t_1})^*\circ D\circ H^* \\
&=&(\text{inj }_{|t_1})^*\circ
\left(((D\otimes\text{id}^*_{\sbb})\circ\mu^*)\times\text{id}^*_{\ca}
\right)\circ F^*\\
&=&  (\text{inj }_{|t_1})^*\circ \left((\mu^*\circ
D)\times\text{id}^*_{\ca}\right)\circ F^*\\
&=&  (\text{inj }_{|t_1})^*\circ \left(\mu^*\times\text{id}^*_{\ca}\right)\circ
F^*\circ X\\
&=& (\text{inj }_{|t_1})^*\circ H^*\circ X .
\end{eqnarray*}
Thus we obtain that $(i)$ implies $(ii)$.\\

Assume now that $(ii)$ is satisfied, i.e., there exists a Lie supergroup
structure on $\sbb$ with multiplication $\mu$ such that 
\begin{equation}
F\circ(\text{id}_{\sbb}\times F)=F\circ(\mu\times\text{id}_{\ca
}).\label{eqaction}
\end{equation}
Since $F$ is a flow for $X$, with initial condition $\phi=\text{id}_{\ca}$, the
LHS of the preceding equality is a flow for $X$ with initial condition $\phi=
F$, (\ref{eqaction}) implies 
\begin{eqnarray}
(\text{inj }_{|t_1})^*\circ(\pa[t_1]+\pa[\tau_1])\circ
(\mu^*\times\text{id}_{\ca}^*)\circ F^*=(\text{inj }_{|t_1})^*\circ
(\mu^*\times\text{id}_{\ca}^*)\circ F^*\circ X.\label{eq34}
\end{eqnarray}
By Equation (\ref{astus2}), the RHS gives for $t_1=0$:
\begin{eqnarray*}
(\text{inj }_{|t_1=0})^*\circ (\mu^*\times\text{id}_{\ca}^*)\circ F^*\circ X &=&
F^*\circ X\\
&=&  (\text{inj }_{\bb}^{\sbb})^*\circ(\pa[t]+\pa[\tau])\circ F^* \\
& &\hspace*{-3cm}+ \tau\cdot\left[ 
F^*\left([X_1,X_0]+\frac{1}{2}[X_1,X_1]\right) +  \pa[\tau]\circ\pa[t]\circ
F^*\right].\\
\end{eqnarray*}
Moreover, we have by direct comparison 
\begin{eqnarray*}
(\pa[t_1]+\pa[\tau_1])\circ \mu^* &=&
(\pa[t_1]+\pa[\tau_1])(\mu^*(t))\cdot(\mu^*\circ\pa[t])\\
&&\hspace*{1.5cm}+(\pa[t_1]+\pa[\tau_1]
)(\mu^*(\tau))\cdot(\mu^*\circ\pa[\tau]).
\end{eqnarray*}
Thus, if $\mu:\sbb\times\sbb\to\sbb$ is given by
\begin{eqnarray*}
\mu^*(t)&=& \widetilde{\mu}(t_1,t_2)+\alpha(t_1,t_2)\tau_1\tau_2 \\
\mu^*(\tau)&=& \beta(t_1,t_2)\tau_1+\gamma(t_1,t_2)\tau_2,
\end{eqnarray*}
and upon using $(\text{inj }_{\{ 0
\}\times\sbb}^{\sbb\times\sbb})^*\circ\mu^*=\text{id}^*_{\sbb}$, we have
\begin{eqnarray*}
(\text{inj }_{\{ 0 \}\times\sbb}^{\sbb\times\sbb})^*(\pa[t_1]+\pa[\tau_1])\circ
\mu^* &=&
\left((\pa[t_1]\widetilde{\mu})(0,t)+\alpha(0,t)\tau\right)\cdot\pa[t]\\
&
&\hspace*{0mm}+\left(\beta(0,t)+(\pa[t_1]\gamma)(0,t)\tau\right)\cdot\pa[\tau].
\end{eqnarray*}
Using again (\ref{astus}), we have
\begin{eqnarray*}
\pa[t]\circ F^* &=& (\text{inj }_{\bb}^{\sbb})^*\circ F^*\circ
X_0+\tau\cdot\pa[\tau]\circ\pa[t]\circ F^*\quad\text{ and }\\
\pa[\tau]\circ F^*&=&(\text{inj }_{\bb}^{\sbb})^*\circ F^*\circ X_1.
\end{eqnarray*}
Then the LHS of (\ref{eq34}) at $t_1=0$ is
\begin{eqnarray*}
&&\hspace*{-2cm}(\text{inj }_{|t_1=0})^*\circ(\pa[t_1]+\pa[\tau_1])\circ
(\mu^*\times\text{id}_{\ca}^*)\circ F^*\\
&&\hspace*{-1cm}=(\pa[t_1]\widetilde{\mu})(0,t)\cdot(\text{inj }_{\bb}^{\sbb})^*\circ F^*\circ
X_0\\
&& + \beta(0,t)\cdot(\text{inj }_{\bb}^{\sbb})^*\circ F^*\circ X_1\\
&& \hspace*{5mm}+ \tau\cdot\big((\pa[t_1]\widetilde{\mu})(0,t)\cdot\pa[\tau]\circ\pa[t]\circ
F^*\\
&& \hspace*{12mm}+
 \alpha(0,t)\cdot  F^*\circ X_0\\
 &&\hspace*{17mm}+ (\pa[t_1]\gamma)(0,t)\cdot F^*\circ X_1\big).
\end{eqnarray*}
Using the obtain identities for its LHS and RHS, the ``$\tau$-part'' of Equation
(\ref{eq34}) at $t_1=0$ gives us: 
\begin{eqnarray*}
&&\hspace*{-.5cm}\tau\cdot\left[  F^*\left([X_1,X_0]+\frac{1}{2}[X_1,X_1]\right) + 
\pa[\tau]\circ\pa[t]\circ F^*\right]\\
&&=
\tau\cdot\big[(\pa[t_1]\widetilde{\mu})(0,t)\cdot\pa[\tau]\circ\pa[t]\circ
F^*+\alpha(0,t)\cdot  F^*\circ X_0
 + (\pa[t_1]\gamma)(0,t)\cdot F^*\circ X_1\big].
\end{eqnarray*}
Since $\widetilde{\mu}(t_1,0)=t_1$, we have $(\pa[t_1]\widetilde{\mu})(0,0)=1$
and therefore the preceding equation evaluated at $t=0$ yields
$$
[X_1,X_0]+\frac{1}{2}[X_1,X_1]=(\pa[t_1]\gamma)(0,0)\cdot X_1+ \alpha(0,0)\cdot
X_0
$$
finishing the proof that $(ii)$ implies $(iii)$.
 \hfill$\square$

\section{Examples and applications}

\begin{noindlist}
 \item \label{Ex1}If $X=X_0$ is an even vector field, the fact that it
integrates to a (local) action of
$\bb=\sbb[1|0]$ is almost folkloristic. The relatively recent proof of $\cite{dumi}$ - in the case of compact
supermanifolds - is close to our approach. A non-trivial (local) action of $\sbb[1|0]$ can obviously be extended to a
(local) action
of $(\sbb,\mu_{a,b})$ if and only if $a=0$. Of course, the ensueing action of $\sbb$ will not even be almost-effective,
since the positive-dimensional sub Lie supergroup $\sbb[0|1]$ acts trivially.

 \item  Our preferred example of an even vector field gives rise to the exponential map on Lie supergroups.

Let us first recall that an even vector field $X$ on a supermanifold $\ca$ corresponds to a section $\sigma_X$ of the
tangent bundle
$T\ca\to\ca$ (see, e.g., Sections 7 and 8 of $\cite{Schmitt}$ for a construction of $T\ca$ and a proof of this 
statement, and compare also the remark after Thm. 2.19 in $\cite{SG-TW}$). Given an auxiliary
supermanifold $\ca[S]$ and a morphism $\psi:\ca[S]\to\ca$, one calls for $i\in\{0,1\}$
\[
 \text{Der}_\psi(\om(M),\om[S](S))_i:=\hspace{7.5cm}\]
\[\{D:\om(M)\to\om[S](S)| D\text{ is }\bb\text{-linear and }\forall
f,g\in\om(M)\text{ homogeneous,}
\]
\[
 D(f\cdot g)=D(f)\cdot \psi^*(g)+\sig[i\cdot|f|]\psi^*(f)\cdot D(g)\},\hspace{3.5cm}
\]
the ``space of derivations of parity $i$ along $\psi$''. In category-theoretical terms the tangent bundle 
$T\ca$ represents then the functor from supermanifolds to sets given by
$ \ca[S]\mapsto\{(\psi,D)|\psi\in\text{Mor}(\ca[S],\ca)\text{ and }D\in\text{Der}_{\psi}(\om(M),\om[S](S))_0\}$
(compare, e.g., Section 3 of \cite{HKST}).\\

Let now $\mathcal{G}=(G,\om[G])$ be a Lie supergroup with multiplication $\mu{=}\mu^{\cg}$ and neutral element $e$. 
We define $X$ in $\text{Der}\left(  \ca[O]_{\ca[G]\times T_e\ca[G] }( G\times T_eG)\right)$ to be the
even vector field on
$ \ca[G]\times T_e\ca[G]$ corresponding to the following section $\sigma_X$ of
$T\ca[G]\times T(T_e\ca[G]) \cong T(\ca[G]\times T_e\ca[G]) \to \ca[G]\times T_e\ca[G]$. We denote the zero-section of
$T\ca[G]\to\ca[G]$ by $\sigma_0$ and the canonical inclusion $T_e\ca[G]\to T\ca[G]$ by $i_e$. Then
$\sigma_X:=(T\mu\circ(\sigma_0\times i_e),0)$, where $T\mu:T\ca[G]\times T\ca[G]\cong T(\ca[G]\times\ca[G])\to T\ca[G]$
is the tangential morphism associated to the multiplication morphism. (For simplicity, we write $0$ for the
zero-section of $T(T_e\ca[G])\to T_e\ca[G]$ here and in the sequel.)\\

We observe that the derivation $X$ corresponding to the section $\sigma_X$ can also be described as follows.
Recall first that we have a canonical map 
$$
\Omega_{\ca[G]}^1(G)\overset{\chi}{\longrightarrow}\ca[O]_{T\ca[G]}(TG),
$$ 
since sections in $\Omega_{\ca[G]}^1(G)$ can be canonically identified with superfunctions on $T\ca[G]$ that
are ``linear in the
fibre directions". Using the inclusion map $i_e:T_e\ca[G]\to T\ca[G]$, the composition $i^*_e\circ\chi$ is a
$\om[G](G)$-linear morphism form $\Omega_{\ca[G]}^1(G)$ to
$\ca[O]_{T_e\ca[G]}(T_eG)$. Here the latter vector space is equipped with the structure of a $\om[G](G)$-module via the
morphism $\pi_e=\pi\circ i_e:T_e\ca[G]\to\ca[G]$ where $\pi:T\ca[G]\to\ca[G]$ is the canonical projection.
Observe
that the
space
$\text{Hom}_{\om[G](G)}\left(
\Omega_{\ca[G]}^1(G),\ca[O]_{T_e\ca[G]}(T_eG) \right)$ is isomorphic to
$\text{Der}_{\pi_e}\left(\om[G](G),\ca[O]_{T_e\ca[G]}(T_eG) \right)$, the space of derivations along $\pi_e$. We denote
by $Y$ the derivation along $\pi_e$ corresponding to the morphism $i^*_e\circ\chi$. The derivation $X$ defined above
via the section $\sigma_X$ is then described as the even vector field in  $\text{Der}\left(  \ca[O]_{\ca[G]\times
T_e\ca[G] }(G\times T_eG
)\right)$  such that 

$$
X(f)=0,\quad \forall f\in \ca[O]_{T_e\ca[G]}(T_eG)\subset \ca[O]_{\ca[G]\times T_e\ca[G]}(G\times T_eG),\text{ and}
$$
$$
X(f)=\left( \text{id}_{\ca[G]}^*\otimes Y\right)\circ\mu^*(f),\quad \forall f\in \ca[O]_{\ca[G]}(G)\subset
\ca[O]_{\ca[G]\times T_e\ca[G]}(G\times T_eG).\vspace{\baselineskip}
$$
Let us  now recall that for $\ca[S]$ an arbitrary supermanifold, and $\phi:\ca\to\can$ a morphism between supermanifolds,
we
have an induced map $\phi(\ca[S]):\ca(\ca[S])=\text{Mor}(\ca[S],\ca)\to\text{Mor}(\ca[S],\can)=\can(\ca[S])$,
$\phi(\ca[S])(\psi):=\phi\circ\psi$. Given a finite-dimensional Lie superalgebra $\mathfrak{g}$ or a Lie supergroup
$\ca[G]$, one easily checks that for all $k\geq 0$, $\mathfrak{g}(\sbb[0|k])$ resp. $\ca[G](\sbb[0|k])$ is a
finite-dimensional classical (i.e. even) Lie algebra resp. Lie group. Furthermore,
 $T_e(\cg(\sk))$ is canonically isomorphic to $(T_e\cg)(\sk)$, where
the first $e$ is the obvious constant morphism from $\sk$ to $\cg$ and the second $e$ denotes the neutral element of
$\cg$. (Compare, e.g., \cite{Sachse} for more information on the superpoint approach to Lie supergroups.)\\

\begin{lem}\label{Lem_exp} Let $\ca[G]$ be a Lie supergroup with multiplication $\mu^{\cg}$,
and the vector field $X$ as above. Then
 \begin{enumerate}[(i)]
  \item \hspace{1mm}the induced vector field $\widetilde{X}$ on the underlying manifold $G\times T_eG$ is given as
$$\widetilde{X}_{(g,\xi)}=(\xi^{L}(g),0)\quad\forall(g,\xi)\in G\times T_eG,$$
\item \hspace{1mm}the (even) vector fields $\widetilde{X}$ and $X$ are complete.
\end{enumerate}
\end{lem}

\textbf{Proof. }(i) For $k\geq 0$, let $\sigma_X(\sk)$ be the section of $T(\cg(\sk))\times
T(T_e\cg(\sk))\to\cg(\sk)\times T_e\cg(\sk)$ induced by
$\sigma_X$, and let $X^k$ be the corresponding derivation on $\cg(\sk)\times T_e\cg(\sk)$. Since $\cg(\sk)\times
(T_e\cg)(\sk)$ is an ungraded manifold, 
\[
 \begin{array}{rcl}
  \sigma_X(\sk)(g,\xi)&=&(T\mu^{\cg}\circ(\sigma_0\times i_e)\circ(g\times\xi),0)\\
		      &=&(T\mu^{\cg}\circ(0_g\times \xi),0)\\
		      &=&\left((T\mu^{\cg(\sk)})_{(g,e)}(0, \xi),0\right)\\
		      &=&\left((Tl^{\cg(\sk)}_g)_{e}(\xi),0\right),
 \end{array}
\]
where $\mu^{\cg(\sk)}$ is the multiplication on $\cg(\sk)$ and $l^{\cg(\sk)}_g$ is the left-multiplication by the
element $g$ of the group $\cg(\sk)$. We conclude that $\sigma_X(\sk)(g,\xi)$ (or equivalently $X^k_{(g,\xi)}$)
corresponds to
$(\xi^L(g),0)$, where $\xi^L$ is the unique left-invariant vector field on $\cg(\sk)$ such that its value in $e$ is
$\xi$. (Observe that $\cg(\sk)$ is a classical Lie group and not only a group object in the category of supermanifolds,
allowing us to argue ``point-wise''.)\\
(ii) The flows of $X^k$ are simply given by $F^{X^k}:\bb\times\cg(\sk)\times(T_e\cg)(\sk)\to\cg(\sk)\times
(T_e\cg)(\sk)$, $(t,g,\xi)\mapsto(g\cdot\text{exp}^{\cg(\sk)}(t\xi),\xi)$. All fields $X^k$ are thus complete, in
particular this holds for $\widetilde{X}=X^0$, the induced vector field on $G=\cg(\sbb[0|0])$. By Theorem
\ref{theoflowm} the flow $F^X:\bb\times\cg\times T_e\cg\to\cg\times T_e\cg$ is then global as well, i.e. $X$ is
complete.\hfill$\square$\\

\begin{de}Let $\cat[G]$ be a Lie supergroup with multiplication $\mu$ and neutral element $e$, and with the even vector
field $X$ and its flow morphism $F=F^X$ as above. Then the ``exponential morphism of $\ca[G]$'' is given by
$exp^{\:\ca[G]}=\text{proj}_1\circ F \circ\text{inj}^{\;\;\;\bb\times\ca[G]\times T_e\ca[G]
}_{\{1\}\times \{e\}\times T_e\ca[G]}:T_e\ca[G] \to \ca[G],
$
where $\text{proj}_1:\ca[G]\times T_e\ca[G]\to\ca[G]$ is the projection on the first factor. Diagrammatically, one has
\[\begin{diagram}
 \node{\hspace{-20mm}\bb\times\ca[G] \times T_e\ca[G]\supset \{1\}\times \{e\} \times  T_e\ca[G]
}\arrow{e,t}{F}\node{\ca[G]\times T_e\ca[G]}\arrow{s,r}{\text{proj}_{1}}\\
\node{T_e\ca[G]}\arrow{n,l}{\cong}\arrow{e,b}{exp^{\ca[G]}}\node{\ca[G].}
\end{diagram}\]
\end{de}

\begin{theo}
 The exponential morphism $exp^{\ca[G]}:T_e\ca[G]\to\ca[G]$ for a Lie supergroup $\ca[G]$
fulfills
and is uniquely determined by the following condition: for all $k\geq 0$,
$exp^{\ca[G]}(\sbb[0|k]):T_e\ca[G](\sbb[0|k])\to\ca[G](\sbb[0|k])$ is the exponential map
$exp^{\ca[G](\sbb[0|k])}$
of the finite-dimensional, ungraded Lie group $\ca[G](\sbb[0|k])$.
\end{theo}
\textbf{Proof. }Using the notations of Lemma \ref{Lem_exp}, a straightforward calculation shows that the flow $F^{X^k}$
of $X^k$ on $(\cg \times T_e\cg)(\sk)$ is given as
follows $(t,(g,\xi))\mapsto F^X\circ\text{inj}_t\circ(g\times \xi)$ and, notably, we have
$F^{X^k}\circ\text{inj}_{(1,e)}(\xi)=F^X\circ\text{inj}_{1}\circ(e\times\xi)=F^X\circ\text{inj}_{(1,e)}\circ \xi$. Hence
\[
 \begin{array}{rcl}
exp^{\cg(\sk)}(\xi)&=&\text{proj}_1\circ F^{X^k}\circ\text{inj}_{(1,e)}(\xi)\\
		     &=&\text{proj}_1\circ F^X\circ\text{inj}_{(1,e)}\circ\xi\\
		     &=&\exp^{\cg}\circ\:\xi\\
		     &=&(\exp^{\cg})(\sk)(\xi).
 \end{array}
\]
On the other hand, it is clear
that the subcategory of superpoints with objects $\{\sbb[0|k]|k\geq0\}$ generates the
category of supermanifolds in the following sense: given two different morphisms $\phi_1,\phi_2:\ca\to\can$ between
supermanifolds, there exists a $k\geq0$ and a morphism $\psi:\sbb[0|k]\to\ca$ such that
$\phi_1\circ\psi\neq\phi_2\circ\psi$. Thus it follows that the family $\{exp^{\ca[G]}(\sbb[0|k])|k\geq0\}$ uniquely
fixes
$exp^{\ca[G]}$.\hfill$\square$

\item Part of our interest in the integration of supervector fields stemmed from the construction of a geodesic
flow in  \cite{SG-TW} . Given a homogeneous (i.e., even or odd) Riemannian metric on a supermanifold $\ca$, the associated geodesic flow is, in fact, defined as the flow of an appropriate Hamiltonian vector field on its 
(co-)tangent bundle.

 \item \label{Ex4}An odd vector field $X_1$ on a supermanifold $\ca$ is called ``homological'' if $X_1\circ
X_1=\frac{1}{2}[X_1,X_1]=0$.
Its flow is given by the following $\sbb[0|1]$-action $\Phi:\sbb[0|1]\times\ca\to\ca$, $\Phi^{*}(f)=f+\tau\cdot X_1(f),
\forall f\in\om(M)$. This action can, of course, be extended to a (not almost-effective) $(\sbb[1|1],\mu_{a,b})$-action
if and
only if $b=0$.\\
Typical examples arise as follows: let $E\to M$ be a vector bundle over a classical manifold and $T$ be an
``appropriate'' $\bb$-linear operator on sections of $\wedge E^*$, then $T$ yields a vector field on $\Pi E:=\left(
M,\Gamma^{\infty}_{\wedge E^*} \right)$, the supermanifold associated to $E\to M$ by the Batchelor
construction. If $E=TM\to M$, we have $\Gamma^{\infty}_{\wedge E^*}=\Omega^\bullet_M(M)$, the sheaf of differential
forms on $M$, with its natural $\bb[Z]/2\bb[Z]$-grading. Taking $T=d$, we get an odd vector field that is obviously
homological. Taking $T=\iota_\xi$, the contraction of differential forms with a vector field $\xi$ on the (here
classical) base manifold $M$, we again get a homological vector field on $\Pi TM$. Since $\iota_\xi\circ\iota_\eta+
\iota_\eta\circ\iota_\xi=0$, the vector space of all vector fields on $M$ is realized as a commutative, purely odd sub
Lie superalgebra of all vector fields on $\Pi TM$. More generally, a section $s$ of $E\to M$ always gives rise to a
contraction $\iota_s:\Gamma^{\infty}_{\wedge
E^*}(M)\to\Gamma^{\infty}_{\wedge E^*}(M)$ that is an odd derivation (i.e. an anti-derivation of degree -1 in more
classical language).
Furthermore, given two sections $s$ and $t$ of $E$, the associated odd vector fields commute. In the article
\cite{ACDS} this construction is studied in the special case that $E$ is the spinor bundle over a
classical
spin manifold $M$.

 \item If $G$ is a Lie group acting on a classical manifold $M$, then the action can of course be lifted to an action on
the total space of the tangent and the cotangent bundle of $M$. The induced vector fields on $\Pi TM$ are even and
$\xi$ in $\mathfrak{g}=\text{Lie}(G)$, the Lie algebra of $G$, acts on $\ca[O]_{\Pi TM}$ by $\ca[L]_\xi$, the Lie
derivative with
respect to the fundamental vector field on $M$ associated to $\xi$. Putting together these
fields and the contractions constructed in Example (4.\ref{Ex4}), we get a Lie superalgebra with underlying vector
space
$\mathfrak{g}\oplus\mathfrak{g}$, the first resp. second summand being the even resp. odd part. The commutators in
$\mathfrak{g}\oplus\mathfrak{g}$ are given as follows: $[\ca[L]_\xi,\ca[L]_\eta]=\ca[L]_{[\xi,\eta]}$,
$[\iota_\xi,\iota_\eta]=0$
and $[\ca[L]_\xi,\iota_\eta]=\iota_{[\xi,\eta]}$ for all $\xi,\eta\in\mathfrak{g}$. In fact, the above can be
interpreted as
an action of the Lie supergroup $\Pi TG$ on $\Pi TM$. \\
The Lie algebra $\mathfrak{g}\oplus\mathfrak{g}$ can be extended by a one-dimensional odd direction generated
by the
exterior derivative $d$. The extended algebra $\mathfrak{g}\oplus(\mathfrak{g}\oplus\bb\cdot d)$, still a sub Lie
superalgebra of $\ca[T]_{\Pi TM}(M)$, has the following additional commutators:
$$[\ca[L]_\xi,d]=0\text{ and }[d,\iota_\xi]=\ca[L]_\xi\text{ for all }\xi\in \mathfrak{g}.$$
 \item If $\ca=\sbb$ with coordinates $(x,\xi)$, the vector field $X_1=\pa[\xi]+\xi\pa[x]$ is obviously odd and
non-homological since $X_1\circ X_1=\pa[x]$. Direct inspection shows that the map
$\Phi:\sbb[0|1]\times\ca\to\ca,\Phi^*(f)=f+\tau\cdot X_1(f)$ (compare Example (4.\ref{Ex1})) does not fulfill
$\pa[\tau]\circ\Phi^*=\Phi^*\circ X_1$. Nevertheless, the trivial extension of $\Phi$ to a morphism
$F:\sbb[1|1]\times\ca\to\ca$ is the flow of $X_1$ in the sense of Theorem \ref{theoflowm}, fulfilling
the initial condition $\phi=\text{id}_{\ca}$.  
We underline that this map is not an action of $\sbb$. Upon extending $X_1$ to $X:=X_0+X_1$ with 
$X_0:=\frac{1}{2}[X_1,X_1]=X_1\circ X_1$, we
obtain by Theorem \ref{theoaction} an action of $(\sbb,\mu_{1,0})$ as the flow map of $X$. Let us observe that the
above vector field $X_1$ (and not $X$) is the prototype of what is called a ``supersymmetry'' in the physics
literature (compare, e.g., \cite{Witten} and the other relevant texts in these volumes.). More recently, the
associated Lie supergroup structure on $\sbb$ (and an analogous structure on $\sbb[2|1]$) were introduced 
by S. Stolz and P. Teichner into their program to geometrize the cocycles of elliptic cohomology (compare \cite{HST} and \cite{ST}).\\

 Obviously, one can generalize this construction to $\sbb[m|n]$ ($m,n\geq 1$) with co\-ordinates
$(x_1,\ldots,x_m,\xi_1,\ldots,\xi_n)$ by setting for $1\leq k\leq m$, $1\leq \alpha\leq n$
$$ D_{\alpha,k}:=\pa[\xi_{\alpha}]+\xi_\alpha\cdot\pa[x_{k}].$$
We then have $[D_{\alpha,k},D_{\beta,l}]=\delta_{\alpha,\beta}\cdot(\pa[x_k]+\pa[x_l])$ and
$[D_{\alpha,k},\pa[x_l]]=0$.
Taking $X_1=D_{\alpha,k}$ and $X_0=\pa[x_k]$ we reproduce a copy of the preceding situation.
\item The vector field $X=X_0+X_1$ on $\ca=\sbb$ with $X_0=\pa[x]+\xi\cdot\pa[\xi]$ and $X_1=\pa[\xi]+\xi\cdot\pa[x]$,
already mentioned in the introduction, is a very simple example of an inhomogeneous vector field not generating any
local $\sbb$-action, since, e.g., condition (iii) in Theorem \ref{theoaction} is violated. Thus integration of $X$
is only possible in the sense of Theorem \ref{theoflowm}, i.e. upon using the evaluation map. Let us observe that the
sub Lie superalgebra $\mathfrak{g}$ of $\tm(M)$ generated by $X$, i.e. by $X_0$ and $X_1$ since sub Lie
superalgebras are by
definition graded sub vector spaces, is four-dimensional with two even generators $Z,W$ and two odd generators
$D,Q$ such that: $Z$ is central, $[W,D]=Q$, $[W,Q]=D$, $D^2=Q^2=Z$ and $[D,Q]=0$. (This amounts in physical
interpretation to the presence of two commuting supersymmetries $D$ and $Q$, generating the same supersymmetric
Hamiltonian $Z$ plus an even symmetry commuting with the Hamiltonian and exchanging the supersymmetries $D$ and
$Q$.)
\end{noindlist}
\section{Flow of a holomorphic vector field on a holomorphic supermanifold}
In this section we extend our results to the holomorphic case. We will always denote the canonical coordinates on
$\bb[C]^{1|1}$ by $z$ and $\zeta$ and write $\pa[z]$ resp. $\pa[\zeta]$ for $\frac{\partial }{\partial z}$ resp.
$\frac{\partial }{\partial \zeta}$.
\begin{de}\label{DefC}
 Let $\ca=(M,\om)$ be a holomorphic supermanifold and $X$ a holomorphic vector field on $\ca$ and $\ca[S]$ a
supermanifold with a morphism $\phi\in\text{Mor}(\ca[S],\ca)$ and $z_0$ in $\bb[C]$.
\begin{enumerate}[(1)]
 \item A ``flow for $X$ (with initial condition $\phi$ and with respect to $z_0$)'' is an open sub supermanifold
$\ca[V]\subset\bb[C]^{1|1}\times\ca[S]$, such that $\{z_0\}\times S\subset V$ ($S$ and $V$ the bodies of $\ca[S]$ and
$\ca[V]$) and such that for all $s$ in $S$, $(\bb[C]\times\{s\})\cap V$ is connected, together with a morphism of
holomorphic supermanifolds $F:\ca[V]\to\ca$ such that
 \begin{eqnarray*}
(\text{inj }_{\bbc}^{\sbbc})^*\circ(\pa[z]+\pa[\zeta])\circ F^* &=&  (\text{inj
}_{\bbc}^{\sbbc})^*\circ F^*\circ X\text{ and }\label{eq2110c}\\
F\circ\text{inj }_{\{z_0\}\times \ca[S]}^{\ca[V]} &=& \phi.\label{eq211c}
\end{eqnarray*}
Sometimes we call the supermanifold $\ca[V]$ (or abusively its body $V$) the ``flow domain'' of $X$.
\item A flow domain $\ca[V]$ of a flow $(\ca[V],F)$ for $X$ is called ``fibrewise 1-connected (relative to the
projection $\ca[V]\to\ca[S]$)'' (or ``fibrewise
1-connected over $\ca[S]$'') if for all $s$ in $S$, $(\bbc\times\{s\})\cap V$ is connected and simply connected.
\end{enumerate}
\end{de}
\textbf{Remark.} We avoid the term ``complex supermanifold'' here, since it is often used to describe supermanifolds
that are, as ringed spaces, locally isomorphic to open sets $D\subset\bb^k$ with structure sheaf
$\ca[C]^\infty_D\otimes_{\bb}\wedge\bbc^l$. ``Holomorphic supermanifolds'' are of course locally isomorphic
to open sets $D\subset\bbc^{k}$ with structure sheaf $\ca[O]_{D}\otimes_{\bbc}\wedge \bbc^l$, where $\ca[O]_{D}$
denotes the sheaf of holomorphic functions on $D$.\\

 Let us first give the holomorphic analogue of Lemma
\ref{theolocal}.
\begin{lem}\label{theolocalC}
 Let $\ca[U]\subset\sbbc[m|n]$ and $\ca[W]\subset\sbbc[p|q]$ be superdomains, $X$ a holomorphic vector field on
$\ca[W]$ (not necessarily homogeneous), $\phi$ in $\text{Mor}(\ca[U],\ca[W])$ and $z_0$ in $\bbc$. Then
\begin{enumerate}[(i)]
 \item it exists a holomorphic flow $(V,\widetilde{F})$ for the reduced holomorphic vector field $\widetilde{X}$ on $U$
with initial condition $\widetilde{\phi}$ and with respect to $z_0$ such that the flow domain $V\subset\bbc\times U$ is
fibrewise 1-connected
over $U$. Furthermore on every flow domain in the sense of Definition \ref{DefC} the holomorphic flow is unique.
\item Let now $(V,\widetilde{F})$ be a fibrewise 1-connected flow domain for $\widetilde{X}$ over $U$. Then there
exists
a unique holomorphic flow $F:\ca[V]\to\ca[W]$ for $X$, with $\ca[V]$ the open sub supermanifold of $\sbbc\times\ca[U]$
with body equal to $V$.
\end{enumerate}
\end{lem}
\textbf{Remark.} The example of the holomorphic vector field $X=(w^2+w^3\xi_1\xi_2)\frac{\partial}{\partial w}$ on
$\ca[W]=\sbbc[1|2]$ with coordinates $(w,\xi_1,\xi_2)$ shows that the condition of fibrewise 1-connectivity of $V$
is not only a technical assumption to our proof. The underlying vector field $\widetilde{X}=w^2\frac{\partial}{\partial
w}$ on $\bbc$, with initial condition $\widetilde{\phi}=\text{id}:\bbc\to\bbc$ with respect to $z_0=0$, can be
integrated to the
flow $\widetilde{F}:V=\bbc^2\backslash\{z\cdot w=1\}\to\bbc$, $\widetilde{F}(z,w)=\frac{1}{1/w-z}$ for $w\neq 0$ and 
$\widetilde{F}(z,0)=0$. Obviously, for $w\neq0$, $(\bbc\times\{w\})\cap V$ is connected, but not simply connected.
Direct inspection now shows that the flow $F$ of $X$ with initial condition $\phi=\text{id}$ and with respect to $z_0=0$
cannot be defined on the whole of
$\ca[V]=(V,\ca[O]_{\sbbc\times\sbbc[1|2]}|_{V})$. \\

\textbf{Proof of Lemma \ref{theolocalC}.} (i) The existence (and the stated unicity property) of a flow
$(\check{V},\widetilde{F})$ for $\widetilde{X}$, with $\{z_0\}\times U\subset \check{V}\subset \bbc\times U$ fulfilling
the initial condition $\widetilde{\phi}$ with respect to $z_0\in\bbc$ is of course a classical application
of the existence of solutions of holomorphic ordinary differential equations (see, e.g., \cite{IY}).
Upon reducing the size of $\check{V}$ we always find flow domains that are fibrewise 1-connected.\\

(ii) The induction procedure of the proof of Lemma \ref{theolocal} can be applied here upon recalling the following
standard facts from the theory of holomorphic linear ordinary differential equations (compare, e.g., \cite{IY}):\\

Fact 1. \textit{Let $\Omega\subset\bbc$ be open and 1-connected (i.e. connected and simply connected), and
$z_0\in\Omega$. If
$A:\Omega\to\text{Mat}(N\times N,\bbc)$ and $b:\Omega\to\bbc^N$ are holomorphic and $\psi_0\in\bbc^N$, then there
exists a unique holomorphic map $\psi:\Omega\to\bbc^N$ fulfilling 
$$\frac{\partial}{\partial z}\psi(z)=A(z)\psi(z)+b(z)$$
such that $\psi(z_0)=\psi_0$.
}\\

Fact 2. \textit{Let $\Omega$ and $z_0$ be as in Fact 1, and let P be a holomorphic manifold (``a parameter space''), and
let $A:\Omega\times P\to\text{Mat}(N\times N,\bbc)$, $b:\Omega\times P\to\bbc^N$, as well as $\psi_0:P\to\bbc^N$ be
holomorphic maps. Then there exists a unique holomorphic map $\psi:\Omega\times P\to \bbc^N$ fulfilling
 $$\frac{\partial}{\partial z}\psi(z,x)=A(z,x)\psi(z,x)+b(z,x)$$
such that $\psi(z_0,x)=\psi_0(x)$, $\forall x\in P$.
}\\

Obviously, to apply these facts in our context, we need the fibrewise 1-connectivity of the ``underlying flow domain''
$V$ for
$\widetilde{X}$.\hfill$\square$\\

Before stating and proving our central result in the holomorphic case, we give the following useful shorthand.
\begin{de}
 Let $\ca[S]$ be a supermanifold, $z_0$ in $\bbc$ and $\ca[N]\subset\sbbc\times\ca[S]$ be an open sub supermanifold
containing $\{z_0\}\times\ca[S]$. Then $\ca[N]^{z_0}$ is defined as the open sub supermanifold of $\cat[N]$ whose body
equals $\displaystyle\underset{s\in S}{\amalg}\left( (\bbc\times\{s\})\cap N \right)^{(z_0,s)}$, where $\left(
(\bbc\times\{s\})\cap N \right)^{(z_0,s)}$ is the connected component of $(\bbc\times\{s\})\cap N$ containing $(z_0,s)$.
\end{de}
\textbf{Remark. } A flow domain $\ca[V]$ in the sense of Definition \ref{DefC}(1) is always open and contains
$\{z_0\}\times\ca[S]$. Furthermore for all $s$ in $ S$ the section
$ (\bbc\times\{s\})\cap V$ is connected. The preceding definition will in fact be useful for discussing intersections
of flow domains in the next theorem.
\begin{theo}
 Let $\ca$ be a holomorphic supermanifold and $X$ a holomorphic vector field on $\ca$, and let $\ca[S]$ be a
holomorphic supermanifold with a holomorphic morphism $\psi:\ca[S]\to\ca$, and $z_0\in\bbc$. Then 
\begin{enumerate}[(i)]
 \item there exists a flow $(V,\widetilde{F})$ for the reduced vector field $\widetilde{X}$ with initial condition
$\widetilde{\phi}$ with respect to $z_0$ such that the flow domain $V\subset\bbc\times S$ is fibrewise 1-connected over
$S$,
\item if $(V,\widetilde{F})$ is as in (i), then there exists a unique flow for $X$ with initial condition $\phi$ with
respect to $z_0$, $F:\ca[V]\to\ca$, where $\ca[V]\subset\sbbc\times\ca[S]$ is the open sub supermanifold with body $V$,
\item if $(\ca[V]_{1},F_1)$ and $(\ca[V]_{2},F_2)$ are two flows for $X$, both with initial condition $\phi$ with
respect to
$z_0$, then $F_1=F_2$ on the flow domain $(\ca[V]_1\cap\ca[V]_2)^{z_0}$,
\item there exists maximal flow domains for $X$ and the germs of their flows coincide on $\{z_0\}\times\ca[S]$.
\end{enumerate}
\end{theo}
\textbf{Proof.} (i) and (ii). It easily follows from Lemma \ref{theolocalC} that $\ca[S]$ can be covered by open
sub supermanifolds $\{\ca[U]^\alpha|\alpha\in A\}$ such that $X_{|\ca[U]^\alpha}$ has a holomorphic flow with initial
condition $\phi_{|\ca[U]_{\alpha}}$ with respect to $z_0$,
$F^\alpha:\ca[V]^\alpha=\Delta_{r_{\alpha}}(z_0)\times\sbbc[0|1]\times\ca[U]^\alpha\to\ca$, where $r_{\alpha}>0$, and
for
$r>0$, $\Delta_{r}(z_0)$ is the open disc of radius $r$ centred in $z_0$. Since $F^\alpha=F^\beta$ on
$\ca[V]^\alpha\cap\ca[V]^\beta$ by the unicity part of Lemma \ref{theolocalC}, we can glue these flows to obtain
$\sbbc\times\ca[S]\supset\dis\ca[V]:=\underset{\alpha\in A}{\cup}\ca[V]^\alpha\overset{F}{\longrightarrow}\ca$, 
a flow for $X$ on $\ca$ with initial condition $\phi$ with respect to $z_0$. Obviously, the ``fibres''
$(\bbc\times\{s\})\cap V$ are 1-connected for all $s$ in $S$, i.e., the flow domain $\ca[V]$ is fibrewise 1-connected
over
$\ca[S]$.\\
Note that if we have a flow $\widetilde{F}$ for the reduced vector field $\widetilde{X}$ on a flow domain $V$ that is
fibrewise 1-connected over $S$, then part (ii) of Lemma \ref{theolocalC} yields a flow for $X$ defined on
$\ca[V]=(V,\ca[O]_{\sbbc\times\ca[S]}|_{V})$.\\

(iii) The body of $(\ca[V]_1\cap\ca[V]_2)^{z_0}$ has as a strong deformation retract the body of $\{z_0\}\times\ca[S]$.
Without loss of generality we can assume that $\ca[S]$ and thus $(\ca[V]_1\cap\ca[V]_2)^{z_0}$ are connected. The local
unicity in Lemma \ref{theolocalC} together with the identity principle for holomorphic morphisms of holomorphic
supermanifolds imply that $F_1=F_2$ on $(\ca[V]_1\cap\ca[V]_2)^{z_0}$.\\

(iv) By Zorn's lemma we get maximal flow domains and by part (iii) the corresponding flows coincide near
$\{z_0\}\times\ca[S]$.\hfill$\square$\\

\textbf{Remarks.} (1) The non-unicity of maximal flow domains for holomorphic vector fields is a well-known phenomenon
already in the ungraded case. A simple example for this is the vector field $X$ on $\bbc^*$ such that 
$X(w)=\frac{1}{w}\frac{\partial}{\partial w}$ for all $w$ in $\bbc^*$.\\
(2) Given the above theorem, the analogues of Lemma \ref{rivf}, Proposition
\ref{LocalAction} and Theorem \ref{theoaction} can now without difficulty be proven to hold for
holomorphic supermanifolds.\\

\textbf{Acknowledgements.} We gratefully acknowledge the support of the SFB/TR 12 Symmetries and Universality in
Mesoscopic Systems program of the Deutsche Forschungsgemeinschaft.

\end{document}